%% file: main.tex
\documentclass[11pt,twoside]{amsart}    

\usepackage{prefix}
 
\title{Coxeter cubes on the graph associahedron and separable elements}
 \author{Harrison Gimenez}
 \address{Texas State University} 
 \email{aws95@txstate.edu}
 \author{SuHo Oh}
 \address{Texas State University} 
 \email{suhooh@txstate.edu}
\date{November 2026}

\begin{document}
\begin{abstract}
A Coxeter cube is a Boolean lattice inside a finite Coxeter group $W$ such that every square is a Coxeter square in the sense of Dyer and Wang. A Coxeter cube of full dimension based at the identity, which we call a based cube, is the same as a partition of the positive roots into inversion sets with one simple root in each block. These partitions were studied by Dimitrov, Gigliotti, Ossip, Paquette and Wehlau under the name of fine decompositions, and by Hohlweg and Pons as partitions of $w_0$. We show that the based cubes of $W$ form a point configuration in the graph associahedron $K(\Gamma)$ of the Coxeter diagram $\Gamma$ of $W$. The vertices of $K(\Gamma)$ are the based cubes whose blocks are parabolic. Every other based cube lies in the interior of an edge or of a cubical face, and we count the based cubes on each edge using the Coxeter arrangement. We also characterize the elements of $W$ that occur as vertices of based cubes. For the cubes at the vertices of $K(\Gamma)$ these are exactly the separable elements of Gaetz and Gao, and for arbitrary based cubes they satisfy a similar condition.
\end{abstract}

\maketitle
\section{Introduction}
Diamonds appear everywhere in the combinatorics of a Coxeter group $W$. Every interval of length two in Bruhat order is a diamond \cite{Deodhar77}, \cite[Lemma~2.7.3]{Bjorner-Brenti05}, a special case of any two reflections generating a dihedral reflection subgroup \cite{Dyer90}. Hence questions about pairs of reflections, pairs of edges in the Bruhat graph, and rank two intervals can be settled inside a dihedral group \cite{Dyer91}. Diamonds also govern the topology of Bruhat intervals \cite{DFLORWW26}.

In the weak order, diamonds only come from commuting generators. Dyer \cite{Dyer11II} and Dyer and Wang \cite{Dyer-Wang19} introduced Coxeter squares. A \emph{Coxeter square} is a quadruple $(w,x,y,z)$ of elements of $W$ with $wx = yz$ and $w(\Phi_x) = \Phi_y$, where $\Phi_x$ denotes the left inversion set of $x$, drawn as a commutative square with these four edge labels. The inversion sets then satisfy $\Phi_{wx} = \Phi_w \sqcup \Phi_y$, and placing the identity at one vertex makes the four vertices a square in the weak order whose sides may have any length. A \emph{Coxeter cube} is a cube with edge labels in $W$ all of whose two-dimensional faces are Coxeter squares.
 
Dyer and Wang ask for a concrete description of these cubes in finite Coxeter groups \cite[\S 5.4, Problem~7]{Dyer-Wang19}. Coxeter cubes were studied by the first author in \cite{Gimenez25}, up to an equivalence called reorientation. The dimension of a Coxeter cube is at most the rank $r$ of $W$, and $r$-cubes exist only when $W$ is finite. A Coxeter cube has labels only on its edges, and our starting point is that placing the identity at the terminal vertex of a Coxeter $r$-cube determines the element at every vertex. We call a Coxeter $r$-cube with the identity at its terminal vertex a \emph{based cube}. A based cube is the same as a partition of $\Phi^+$ into $r$ inversion sets with one simple root in each block. In type $A_n$ the Coxeter $n$-cubes modulo reorientation are in bijection with binary trees with $n+1$ leaves modulo isomorphism \cite{Gimenez25}, and placing the identity at the terminal vertex unwraps the equivalence from both sides. The based cubes are in bijection with the vertices of the associahedron, and hence with planar binary trees with $n+1$ leaves.
 
Based cubes have been studied under other names. Partitions of inversion sets into inversion sets are called decompositions of inversion sets in \cite{DDMRWW17}, inv-decompositions in \cite{Katthan13} and partitions of $w$ in \cite{Hohlweg-Pons26}. The based cubes are the finest such partitions of $\Phi^+$, called fine decompositions by Dimitrov, Gigliotti, Ossip, Paquette and Wehlau \cite{Dimitrov-etal25}, who give a recursive description of them for finite Weyl groups. These decompositions arise in geometric invariant theory, in the study of the Belkale--Kumar product \cite{Belkale-Kumar06} on the cohomology $H^*(G/B,\mathbb{Z})$ of the complete flag variety $G/B$ of a complex semisimple algebraic group $G$ \cite{Dimitrov-Roth09,Dimitrov-Roth17}. In type $A$ they are also related to the Babington--Smith model in algebraic statistics \cite{Katthan13} and to the simplicial faces of the Littlewood--Richardson cone \cite{DDMRWW17}.
 
The main result of this paper is that the based cubes of $W$ form a point configuration in the graph associahedron of the Coxeter diagram of $W$. For a graph $\Gamma$, a \emph{tube} is a connected set of nodes and a \emph{tubing} is a collection of tubes that are pairwise nested or disjoint and non-adjacent, where we always include the whole node set as a tube. The \emph{graph associahedron} $K(\Gamma)$ \cite{Carr-Devadoss06,Postnikov09} is a simple polytope whose faces are indexed by the tubings of $\Gamma$. The vertices are the maximal tubings, and each face is a product of smaller graph associahedra. For the path, cycle and complete graph one obtains the associahedron, the cyclohedron and the permutohedron respectively. When $\Gamma$ is the Coxeter diagram of $W$, Ardila, Reiner and Williams \cite{Ardila-Reiner-Williams06} showed that $K(\Gamma)$ is dual to the positive Bergman complex of the Coxeter arrangement of $W$.
 
We take $\Gamma$ to be the Coxeter diagram of $W$, so that a tube $T$ is a connected subset of $S$ and gives the parabolic subgroup $W_T$. A face of $K(\Gamma)$ is \emph{cubical} if it is a product of edges. Each edge of $K(\Gamma)$ determines a tube $T$ and two nodes $p, q \in T$, and we write $\sigma(T,p,q)$ for the number of atoms of a rank two interval in the intersection lattice of the Coxeter arrangement of $W_T$ given by this data. We have $\sigma(T,p,q) \ge 3$, with equality in type $A$.

\begin{theorem*}[\cref{thm:nestedcubesgamma,thm:cubicaltubing}]
The based cubes of $W$ form a point configuration in $K(\Gamma)$ with the following properties.
\begin{enumerate}
\item The vertices of $K(\Gamma)$ are in bijection with the based cubes with all blocks being parabolic inversion sets, which we call \emph{nested cubes}.
\item Every other based cube lies in the interior of a cubical face. In particular no based cube lies in the interior of a face that is not a product of edges.
\item An edge of $K(\Gamma)$ contains exactly $\sigma(T,p,q) - 3$ based cubes in its interior. They are totally ordered along the edge, and two consecutive ones differ by moving a single set of roots from one block to another.
\item A cubical face is a product of its edges, and the points in its interior are exactly the products of the points in the interiors of its edges.
\end{enumerate}
\end{theorem*}

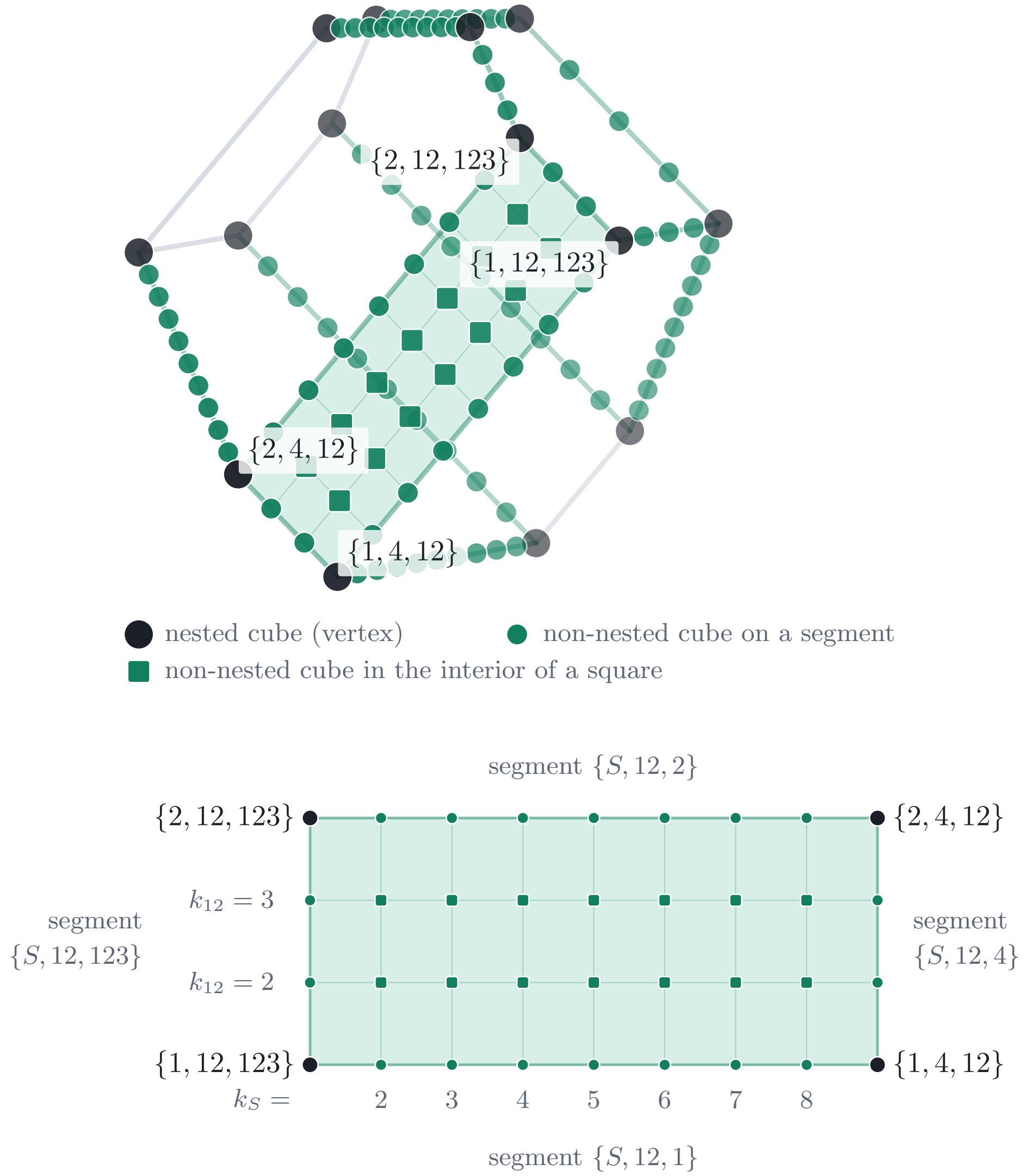
\begin{figure}[ht]
\centering
\input{fig/assoc-H4}
\\[1.5em]
\input{fig/square-H4}
\caption{Top figure shows the $118$ based cubes of $H_4$ as points of the three-dimensional associahedron $K(\Gamma)$, as in \cref{ex:H4}. Edges with $\sigma = 3$ are grey and edges with $\sigma > 3$ are green, with $\sigma - 3$ equally spaced non-nested cubes on each edge. The shaded square is the face $\{S, 12\}$, where the vertices are labelled by their maximal tubings with the tube $S$ omitted. Bottom figure shows the face $\{S, 12\}$ drawn separately. In its interior we have a $7 \times 2$ grid of based cubes parametrized by $(k_S, k_{12})$ as in \cref{thm:cubicaltubing}.}
\label{fig:H4}
\end{figure}

Since the based cubes on an edge are placed at equal distances, the point configuration is canonical. In type $A$ every edge has $\sigma(T,p,q) = 3$, so the configuration consists of the vertices alone and every based cube is nested. Points in the interiors of cubical faces only appear outside type $A$. For type $H_4$, $K(\Gamma)$ is the three-dimensional associahedron and there are $118$ based cubes, $14$ at the vertices, $90$ on the edges and $14$ on one of the three square faces (\cref{fig:H4}). Summing over the cubical faces refines the counts of based cubes obtained in \cite{Dimitrov-etal25} for Weyl groups (\cref{cor:count}).
 
Finally, we characterize the elements of $W$ that occur as vertices of based cubes. For nested cubes these are exactly the separable elements of Gaetz and Gao \cite{Gaetz-Gao20sep} (\cref{thm:nestedvertices}). We introduce \emph{loosely separable} elements, obtained by dropping the parabolic condition, and show that they are exactly the vertices of arbitrary based cubes (\cref{thm:vertices}).
 
The paper is organized as follows. In \cref{sec:prelim} we recall the background on root systems, inversion sets, the Coxeter arrangement, graph associahedra and separable elements. In \cref{sec:cubes} we recall Coxeter squares and Coxeter cubes, define based cubes and characterize them as partitions of $\Phi^+$. In \cref{sec:nested} we study nested cubes, prove their bijection with the vertices of $K(\Gamma)$ and identify their vertices with separable elements. In \cref{sec:nonnested} we treat the remaining based cubes, place them on the cubical faces of $K(\Gamma)$, count them via the Coxeter arrangement and characterize their vertices. Finally \cref{sec:questions} contains further questions.

\section{Preliminaries}
\label{sec:prelim}
 
\subsection{Coxeter groups and root systems}
\label{subsec:roots}
 
Throughout the paper, $W$ is a finite Coxeter group with simple generating set $S$. In other words, $S$ is a finite set and $W$ is the group generated by $S$ subject to the relations
\[
(st)^{m_{st}} = e
\]
for some $m_{st} \in \mathbb{Z}_{>0}$ where $m_{st} = 1$ if and only if $s = t$. We say an expression of $w \in W$ in the simple generators
\[
w = s_1 \cdots s_k
\]
is \emph{reduced} if $w$ cannot be expressed in fewer generators. The length of any reduced expression is unique, so we define the function $\ell : W \to \mathbb{Z}_{\geq 0}$ which maps $w \in W$ to the length of any reduced expression. We call the value $\ell(w)$ the \emph{length} of $w$. We write $r = |S|$ for the rank of $W$ and $w_0$ for the longest element of $W$. The \emph{Coxeter diagram} $\Gamma$ of $W$ is the labelled graph with node set $S$ and edges $(s,t)$ labelled by the value $m_{st}$ under the conventions that we draw no edge if $m_{st} = 2$ and an unlabelled edge if $m_{st} = 3$. For more details, we refer readers to \cite{Bjorner-Brenti05,Humphreys90,Oh-Richmond26}.
 
Each finite Coxeter group $W$ is naturally associated with a hyperplane arrangement through its root system, for which we use the geometric representation of $W$ \cite[\S 5.3]{Humphreys90}. Let $V$ be the real vector space with basis $\Delta := \{\alpha_s \mid s \in S\}$ and the symmetric bilinear form determined by $(\alpha_s,\alpha_t) = -\cos(\pi/m_{st})$. Since $W$ is finite this form is an inner product \cite[Theorem~6.4]{Humphreys90}, and $W$ acts faithfully on $V$ by letting $s \in S$ act as the reflection $s(x) = x - 2(\alpha_s,x)\,\alpha_s$. The elements of $\Delta$ are the \emph{simple roots} and the \emph{root system} of $W$ is
\[
\Phi := \{w\alpha_s \mid w \in W,\ s \in S\},
\]
a finite set of unit vectors with $\Phi = -\Phi$ and $\dim V = r$. Every root is either a nonnegative or a nonpositive linear combination of simple roots \cite[\S 5.4]{Humphreys90}, so we can decompose
\[
\Phi = \Phi^+ \sqcup \Phi^-
\]
into \emph{positive} and \emph{negative} roots. For $\alpha \in \Phi$ let
\[
s_\alpha(x) := x - \frac{2(\alpha,x)}{(\alpha,\alpha)}\,\alpha
\]
be the reflection in the hyperplane orthogonal to $\alpha$. Then $s_{w\alpha} = w s_\alpha w^{-1}$, so $s_\alpha \in W$ for every $\alpha \in \Phi$, and the $s_\alpha$ with $\alpha \in \Phi^+$ are exactly the reflections of $W$ \cite[\S 5.7]{Humphreys90}. The notions we use depend only on the lines spanned by the roots, so each root may be rescaled by a positive constant. In particular, when $W$ is a Weyl group we use a crystallographic root system in the examples.
 
To each $\alpha \in \Phi$, there is the corresponding hyperplane given by
\[
H_\alpha := \{x \in V \mid (\alpha,x) = 0\}.
\]
Note that $H_\alpha = H_{-\alpha}$ and that if $x \in H_\alpha$, then $s_\alpha(x) = x$. The collection of hyperplanes
\[
\A := \{H_\alpha \mid \alpha \in \Phi^+\}
\]
is called the \emph{Coxeter arrangement} of $W$. As each $H_\alpha$ contains the origin in $V$, $\A$ is a central hyperplane arrangement. We follow \cite[Chapter~2]{Orlik-Terao92} and \cite[\S 3.11]{Stanley12} for terminology on hyperplane arrangements. A \emph{flat} of $\A$ is an intersection of hyperplanes of $\A$, and the flats ordered by reverse inclusion form the \emph{intersection lattice} $\LL(\A)$. The \emph{chambers} of $\A$ are the connected components of $V \setminus \bigcup_\beta H_\beta$. For $\beta \in \Phi^+$ let $H_\beta^+ = \{\lambda \in V : (\beta,\lambda) > 0\}$ and $H_\beta^- = -H_\beta^+$. We let
\[
C = \{\lambda \in V : (\alpha_s,\lambda) > 0 \text{ for all } s \in S\} = \bigcap_{\beta \in \Phi^+} H_\beta^+
\]
be the \emph{fundamental chamber}.
 
The \emph{support} of a root $\beta = \sum_{s \in S} c_s \alpha_s$ is $\supp(\beta) = \{s \in S : c_s \neq 0\}$. For $I \subseteq S$ let $W_I$ be the standard parabolic subgroup generated by $I$, and put $\Delta_I := \{\alpha_s : s \in I\}$. Then $\Phi_I := \Phi \cap \spa(\Delta_I)$ is a root system of $W_I$ with simple system $\Delta_I$ \cite[\S 1.10]{Humphreys90}, and has positive roots
\[
\Phi_I^+ = \Phi_I \cap \Phi^+ = \{\beta \in \Phi^+ : \supp(\beta) \subseteq I\}.
\]
The support of every root is connected in $\Gamma$, and every connected subset $T \subseteq S$ is the support of some root \cite[Lemma~5.3]{Ardila-Reiner-Williams06}.
 
\subsection{Inversion sets}
\label{subsec:biclosed}
 
For each element $w \in W$, we can take a subset of positive roots corresponding to the inversions of $w$. We define the \emph{inversion set} of $w$ as the set of positive roots
\[
\Phi_w := \Phi^+ \cap w(\Phi^-) = \{\beta \in \Phi^+ \mid w^{-1}(\beta) \in \Phi^-\}.
\]
 
We have $|\Phi_w| = \ell(w)$ \cite[\S 1.6]{Humphreys90}, $\Phi_e = \emptyset$, $\Phi_s = \{\alpha_s\}$ for $s \in S$, and $\Phi_{w_0} = \Phi^+$. Since $w_0(\Phi^-) = \Phi^+$ we have $\Phi_{ww_0} = \Phi^+ \cap w(\Phi^+)$ and the following equation
\begin{equation}\label{eqn:complement}
\Phi_{ww_0} = \Phi^+ \setminus \Phi_w \qquad \text{for all } w \in W.
\end{equation}
A simple root $\alpha_s$ lies in $\Phi_w$ if and only if $s$ is a left descent of $w$ ($\ell(sw) < \ell(w)$) \cite[\S 1.6]{Humphreys90}. In particular $\Phi_w$ contains a simple root whenever $w \neq e$.
 
Inversion sets are described geometrically by the chambers of $\A$.
 
\begin{proposition}[{\cite[\S 1.12]{Humphreys90}, \cite[Proposition~2.11]{Hohlweg-Labbe16}}]
\label{prop:chambers}
Let $\A$ be the Coxeter arrangement of a finite Coxeter group $W$.
\begin{enumerate}
\item The chambers of $\A$ are the cones $vC$, $v \in W$, and these are pairwise distinct.
\item For $v \in W$ and $\lambda \in vC$,
\[
\Phi_v = \{\beta \in \Phi^+ : (\beta,\lambda) < 0\} = \{\beta \in \Phi^+ : vC \subseteq H_\beta^-\},
\]
the set of positive roots whose hyperplanes separate $C$ from $vC$.
\item A subset of $\Phi^+$ is an inversion set if and only if it is of the form $\{\beta \in \Phi^+ : (\beta,\lambda) < 0\}$ for some $\lambda \in V$ lying on no hyperplane of $\A$.
\end{enumerate}
\end{proposition}

The \emph{(right) weak order} on $W$ is defined by $u \leq_R v$ if $v = ux$ for some $x \in W$ with $\ell(v) = \ell(u) + \ell(x)$. By \cite[Proposition~3.1.3]{Bjorner-Brenti05},
\begin{equation}\label{eqn:weakorder}
u \leq_R v \iff \Phi_u \subseteq \Phi_v ,
\end{equation}
so $w \mapsto \Phi_w$ is an isomorphism from the weak order onto the collection of inversion sets ordered by inclusion, and in particular $w$ is determined by $\Phi_w$. We will repeatedly use the following special case of the cocycle identity for inversion sets \cite[Chapter~1, Exercise~12]{Bjorner-Brenti05}, \cite[\S 2]{Dyer90}.
 
\begin{lemma}\label{lem:cocycle}
If $u, v \in W$ satisfy $\ell(uv) = \ell(u) + \ell(v)$, then $u(\Phi_v) \subseteq \Phi^+$ and
\begin{equation}\label{eqn:cocycleroots}
\Phi_{uv} = \Phi_u \sqcup u(\Phi_v).
\end{equation}
\end{lemma}
 
For $I \subseteq S$ let $w_0(I)$ denote the longest element of $W_I$. Every $w \in W_I$ permutes $\Phi^+ \setminus \Phi_I^+$ \cite[\S 1.10]{Humphreys90}. Consequently, for $w \in W_I$ the inversion set $\Phi_w$ lies in $\Phi_I^+$ and coincides with the inversion set of $w$ computed in the Coxeter group $W_I$ with root system $\Phi_I$. Applying \eqref{eqn:complement} inside $W_I$ on $\Phi_{w_0(I)} = \Phi_I^+$ tells us that
\begin{equation}\label{eqn:paraboliccomplement}
\Phi_{w w_0(I)} = \Phi_I^+ \setminus \Phi_w \qquad \text{for all } w \in W_I.
\end{equation}
 
For $I \subseteq S$ and $w \in W$ there is a unique element $w|_I \in W_I$ with
\begin{equation}\label{eqn:restriction}
\Phi_{w|_I} = \Phi_w \cap \Phi_I^+ ,
\end{equation}
coming from the parabolic decomposition $w = w_I \cdot {}^I w$ \cite[Proposition~2.4.4]{Bjorner-Brenti05}. We call $w|_I$ the \emph{restriction} of $w$ to $I$.

The following lemma is proved for finite Weyl groups in \cite[Proposition~2.12]{Dimitrov-etal25} and for arbitrary Coxeter systems in \cite[Lemma~2.8]{Hohlweg-Pons26}. 
 
\begin{lemma}[{\cite[Lemma~2.8]{Hohlweg-Pons26}, \cite[Proposition~2.12]{Dimitrov-etal25}}]
\label{lem:subunion}
Let $\Phi_w = B_1 \sqcup \cdots \sqcup B_k$ be a partition of an inversion set into inversion sets. Then $\bigcup_{i \in I} B_i$ is an inversion set for every $I \subseteq \{1, \dots, k\}$.
\end{lemma}

 \subsection{Rank two flats}
\label{subsec:ranktwo}

In this subsection we study the two-dimensional flats of $\A$ coming from the simple roots, which are the main tool of \cref{sec:nonnested}. For $J \subseteq S$ let
\[
X_J := \bigcap_{s \in J} H_{\alpha_s} = \spa(\Delta_J)^\perp ,
\]
a flat of dimension $r - |J|$. The positive roots $\beta$ vanishing on $X_J$ ($X_J \subseteq H_\beta$) are exactly $\Phi_J^+$ by \cite[\S 1.10]{Humphreys90}, and every other positive root has a constant nonzero sign on each connected component of $X_J \setminus \bigcup_{\beta \notin \Phi_J^+} H_\beta$.

Now let $J = S \setminus \{p,q\}$, so that $X_J$ is a plane, and put $U := \Phi^+ \setminus \Phi_J^+$. For $\lambda \in X_J$ we will write $\lambda_p := (\alpha_p,\lambda)$ and $\lambda_q := (\alpha_q,\lambda)$. These are linear coordinates on $X_J$, and for $\beta = \sum_{s \in S} c_s(\beta)\,\alpha_s \in \Phi^+$ we have
\begin{equation}\label{eqn:planecoords}
(\beta,\lambda) = c_p(\beta)\,\lambda_p + c_q(\beta)\,\lambda_q .
\end{equation}
Hence for $\beta \in U$ the line $H_\beta \cap X_J$ depends only on the ratio
\[
\rho(\beta) := c_q(\beta)/c_p(\beta) \in [0,\infty] ,
\]
and $H_\beta \cap X_J = H_\gamma \cap X_J$ if and only if $\rho(\beta) = \rho(\gamma)$. Let $\rho_1 < \cdots < \rho_\sigma$ be the values of $\rho$ on $U$ and let $L_i := \{\beta \in U : \rho(\beta) = \rho_i\}$, so that $U = L_1 \sqcup \cdots \sqcup L_\sigma$. We have $\rho_1 = 0$ and $\rho_\sigma = \infty$, where $L_1$ is the class of $\alpha_p$ and $L_\sigma$ the class of $\alpha_q$. 

For $\lambda \in X_J$ we define
\[
B_\lambda := \{\beta \in U : (\beta,\lambda) < 0\},
\]
the set of roots of $U$ such that the corresponding hyperplane $H_\beta$ separates the fundamental chamber from $\lambda$ in $\A$. By \eqref{eqn:planecoords}, if $\lambda_p < 0 < \lambda_q$, then $B_\lambda = \{\beta \in U : \rho(\beta) < -\lambda_p/\lambda_q\}$ can be partitioned as $L_1 \sqcup \cdots \sqcup L_k$. Moreover, for any $k$ with $1 \le k \le \sigma - 1$, we can pick $-\lambda_p/\lambda_q$ strictly between $\rho_k$ and $\rho_{k+1}$ to obtain that particular partition. Similarly $\lambda_p > 0 > \lambda_q$ gives the sets $L_{k+1} \sqcup \cdots \sqcup L_\sigma$, while $\lambda_p, \lambda_q > 0$ gives $B_\lambda = \emptyset$ and $\lambda_p, \lambda_q < 0$ gives $B_\lambda = U$.

\begin{lemma}\label{lem:plane}
Let $J = S \setminus \{p,q\}$ and $U = \Phi^+ \setminus \Phi_J^+$.
\begin{enumerate}
\item If $\lambda \in X_J$ lies on no hyperplane $H_\beta$ with $\beta \in U$, then $B_\lambda$ and $B_\lambda \sqcup \Phi_J^+$ are inversion sets.
\item If $B \subseteq U$ is such that $B$ and $B \sqcup \Phi_J^+$ are inversion sets, then $B = B_\lambda$ for some $\lambda \in X_J$ lying on no hyperplane $H_\beta$ with $\beta \in U$.
\end{enumerate}
\end{lemma}

\begin{proof}
For a flat $X$ of $\A$ write $\A_X$ for the set of hyperplanes of $\A$ containing $X$. The point $\lambda$ lies on $H_\beta$ exactly for $\beta \in \Phi_J^+$, so the chambers with $\lambda$ in their closure are separated from $C$ by the hyperplanes $H_\beta$ with $\beta \in B_\lambda$ together with those with $\beta$ in an inversion set of $W_J$. Taking the two chambers corresponding to $\emptyset$ and $\Phi_J^+$ gives (1) by \cref{prop:chambers}(3).

For the second claim, let $uC$ be the chamber that is separated from the fundamental chamber by the hyperplanes $H_\beta$ with $\beta \in B$, and let $vC$ be the chamber corresponding to $B \sqcup \Phi_J^+$. The hyperplanes that separate $uC$ and $vC$ are exactly $\A_{X_J}$. Translating this by $u^{-1}$, the hyperplanes that separate $C$ and $xC$ (we set $x := u^{-1}v$) are exactly the hyperplanes of $\A_Y$ (hyperplanes containing $Y := u^{-1}(X_J)$). Let $K$ denote the set of simple reflections $s$ that are left descents of $x$, which are those with $H_{\alpha_s} \in \A_Y$. Then $Y \subseteq X_K$, so all hyperplanes containing $X_K$ separate $C$ from $xC$. This means that we can find a walk of the chambers (going through codimension one faces, crossing each hyperplane at most once) from $C$ to $xC$ while passing through $w_0(K)C$. Hence $x = w_0(K)y$ with $\ell(x) = \ell(w_0(K)) + \ell(y)$, and \eqref{eqn:cocycleroots} gives $\Phi_x = \Phi_K^+ \sqcup w_0(K)(\Phi_y)$.

Assume for the sake of contradiction that $y \neq e$ and let $s$ be a left descent of $y$, so that $\alpha_s \in \Phi_y$. Then $w_0(K)\alpha_s \in w_0(K)(\Phi_y) \subseteq \Phi_x$ is a positive root. Since $w_0(K)$ sends $\alpha_t$ to a negative root for every $t \in K$, we get $s \notin K$. Moreover $H := w_0(K)H_{\alpha_s}$ is a hyperplane that separates $C$ from $xC$, implying $H \in \A_Y$. As $w_0(K)$ fixes all points of $X_K$ (containing $Y$), it permutes $\A_Y$, so $H_{\alpha_s} = w_0(K)H \in \A_Y$. However this means $s \in K$ and we get a contradiction.

From the above we can conclude that $xC = w_0(K)C$, so we have $\A_Y = \A_{X_K}$. Since two flats that are contained in the same set of hyperplanes have to be equal, we get $Y = X_K$ and $X_J = u(X_K)$.

Pick $\mu \in X_K$ on the positive side of all hyperplanes $H_{\alpha_s}$ for $s \notin K$, and put $\lambda := u\mu \in X_J$. Let $\beta \in U$. Then $H_\beta \notin \A_{X_J}$, so $H_{u^{-1}\beta} \notin \A_{X_K}$ and $u^{-1}\beta \notin \Phi_K$. Hence $u^{-1}\beta$ has a positive coefficient on some $\alpha_s$ with $s \notin K$ if it is a positive root, and a negative one if it is a negative root. As $\mu$ vanishes on $\Delta_K$ and is positive on the other simple roots, $(\beta,\lambda) = (u^{-1}\beta,\mu)$ is nonzero, and it is negative if and only if $u^{-1}\beta$ is a negative root (equivalent to $\beta \in \Phi_u = B$). Hence $\lambda$ lies on no $H_\beta$ with $\beta \in U$ and $B_\lambda = B$.
\end{proof}

\subsection{Graph associahedra}
\label{subsec:assoc}
We follow \cite{Carr-Devadoss06}, \cite{Postnikov09} and \cite{Gaetz-Gao20split} for the terminology on graph associahedra. Let $\Gamma$ be a graph on the node set $S$. A \emph{tube} is a nonempty subset $T \subseteq S$ inducing a connected subgraph, and two tubes are \emph{compatible} if
\begin{itemize}
    \item they are nested, or
    \item they are disjoint and their union is not a tube. That is, no edge of $\Gamma$ joins them.
\end{itemize}
A \emph{tubing} is a set $N$ of pairwise compatible tubes. This is the definition of \cite{Carr-Devadoss06,Ardila-Reiner-Williams06}, except that we also allow the connected components of $\Gamma$ to be tubes (as in \cite{Postnikov09,Gaetz-Gao20split}). Since no edge of $\Gamma$ joins two disjoint compatible tubes, tubings containing all connected components of $\Gamma$ are exactly the nested sets (union of any two or more pairwise disjoint tubes is not a tube) of \cite[Definition 7.3]{Postnikov09}.

 The \emph{graph associahedron} $K(\Gamma)$ is a simple polytope of dimension $|S| - c$ ($c$ is the number of connected components of $\Gamma$) with face poset the poset of tubings containing all connected components of $\Gamma$ (ordered by reverse inclusion). The vertices of $K(\Gamma)$ correspond to the inclusion-wise maximal tubings, which all have $|S|$ elements, and a tubing $N$ corresponds to a face $\mathcal{F}_N$ of dimension $|S| - |N|$. When $\Gamma$ is a path on $n$ vertices, $K(\Gamma)$ is the usual associahedron.

We only consider $\Gamma$ the Coxeter diagram of $W$. For irreducible $W$ this is a tree, and for reducible $W$ the graph associahedron is the product of the graph associahedra of the components. Unless otherwise stated, $\Gamma$ is a tree, so that $S$ is a tube and every maximal tubing contains $S$.

Given a tubing $N$, the \emph{children} of $T \in N$ are the maximal tubes of $N$ properly contained in $T$. The \emph{parent} of $T \neq S$ is the minimal tube of $N$ properly containing $T$. We define the \emph{shell} of $T$ as
\[
\shell_N(T) = T \setminus \bigcup \{T' : T' \text{ a child of } T\},
\]
the set of nodes of $T$ lying in no smaller tube of $N$. We omit $N$ when $N$ is clear from context. The shells of the tubes of $N$ partition $S$, and $N$ is maximal if and only if every shell is a single node. In such cases we use $p(T)$ to denote that unique node of the shell and call it the \emph{pivot}, following \cite{Gaetz-Gao20sep}.

A tubing $N$ with exactly one shell $\shell_N(T) = \{x,y\}$ of size two and all other shells of size one corresponds to an edge of $K(\Gamma)$. Its two endpoints are the maximal tubings $N \cup \{T_x\}$ and $N \cup \{T_y\}$, where $T_x$ is the connected component of $T \setminus \{y\}$ containing $x$ and $T_y$ is the connected component of $T \setminus \{x\}$ containing $y$. Two maximal tubings that are the endpoints of an edge are said to differ by a \emph{flip}. Such tubings differ by exactly one tube. 

For any tubing $N$ and any $T \in N$ with $\shell_N(T) = \{x,y\}$ we say that $N \cup \{T_x\}$ and $N \cup \{T_y\}$ are obtained from $N$ by \emph{completing} the shell of $T$. These are again tubings. In $N \cup \{T_x\}$ the shell of $T$ is $\{y\}$ and the shell of $T_x$ is $\{x\}$ while all other shells are unchanged. Similar for $N \cup \{T_y\}$.

By \cite[Proposition 7.5]{Postnikov09} and its proof, the face $\mathcal{F}_N$ of $K(\Gamma)$ corresponding to a tubing $N$ is a product over $T \in N$ of polytopes of dimension $|\shell_N(T)| - 1$. We call a face (together with its tubing) \emph{cubical} if every shell of the tubing has size at most two. Such a face is combinatorially a cube whose dimension is the number of size two shells. Its vertices are the maximal tubings obtained from $N$ by completing every shell of size two and its edges are obtained by completing all but one of them.

\subsection{Separable elements}
Following \cite{Gaetz-Gao20sep}, write $I_\Phi(w) = \{\beta \in \Phi^+ : w\beta \in \Phi^-\}$, so that $I_\Phi(w) = \Phi_{w^{-1}}$ in our notation. For $\Delta' \subseteq \Delta$ generating the subsystem $\Phi'$, the restriction $w|_{\Phi'}$ is the element of $W(\Phi')$ with $I_{\Phi'}(w|_{\Phi'}) = I_\Phi(w) \cap \Phi'$.

\begin{definition}[{cf.\ \cite[Definition 3.1]{Gaetz-Gao20sep}}]
Let $w \in W(\Phi)$. Then $w$ is \emph{separable} if one of the following holds:
\begin{itemize}
    \item $\Phi$ is of type $A_1$;
    \item $\Phi = \bigoplus \Phi_i$ is reducible and $w|_{\Phi_i}$ is separable for each $i$;
    \item $\Phi$ is irreducible and there exists a \emph{pivot} $\alpha_s \in \Delta$ such that $w|_{\Phi'} \in W(\Phi')$ is separable where $\Phi'$ is generated by $\Delta' = \Delta \setminus \{\alpha_s\}$ and such that either
    \[
    \{\beta \in \Phi^+ : s \in \supp(\beta)\} \subseteq I_\Phi(w), \quad \text{or}
    \]
    \[
    \{\beta \in \Phi^+ : s \in \supp(\beta)\} \cap I_\Phi(w) = \emptyset.
    \]
\end{itemize}
Since this definition depends not only on $w$ but also on $\Phi$, we sometimes say that $(w,\Phi)$ is separable for clarity.
\end{definition}

In the third clause, following \cite{Gaetz-Gao20split,Gaetz-Gao26}, we call $\alpha_s$ a \emph{full pivot} of $w$ in the first case (when $\{\beta \in \Phi^+ : s \in \supp(\beta)\} \subseteq I_\Phi(w)$), and an \emph{empty pivot} in the second case (when this set is disjoint from $I_\Phi(w)$). We also refer to the corresponding node of $\Gamma$ as a pivot.

Gaetz and Gao state the definition for finite Weyl groups, with the set of roots $\beta \in \Phi^+$ such that $\beta - \alpha_s$ is a nonnegative combination of simple roots instead of $\{\beta \in \Phi^+ : s \in \supp(\beta)\}$. For a crystallographic root system the two sets are the same, since the coefficients of the roots are integers. Our definition makes sense verbatim for any finite Coxeter group (which is the scope of our paper) and does not depend on the normalization of the roots. We also read the definition with $\Phi_w$ in place of $I_\Phi(w)$: we call $w$ separable when $w^{-1}$ is separable in the sense above. The two conventions agree, since the set of separable elements is closed under $w \mapsto w^{-1}$: for finite Weyl groups this is \cite[Corollary 2]{Gaetz-Gao20split}, for $H_3$ and $H_4$ it was checked by computer \cite{GO26code}, and for $I_2(m)$ the six separable elements are $e$, $w_0$ and the four elements of length $1$ and $m-1$. The restriction $w|_T$ of \eqref{eqn:restriction} corresponds to the restriction of Gaetz and Gao in the same way: $(w|_T)^{-1} = (w^{-1})|_{\Phi_T}$.

In type $A$ the separable elements are exactly the permutations avoiding $3142$ and $2413$, which were studied in \cite{Shapiro-Stephens91,West95,Bose-Buss-Lubiw98}. When $w$ is separable, the weak order intervals $[e,w]$ and $[w,w_0]$ are rank-symmetric and rank-unimodal, and the product of their rank generating functions is the rank generating function of $W$ \cite[Theorem 3.9]{Gaetz-Gao20sep}. By \cite[Theorem 7]{Gaetz-Gao20split} the separable elements of $W$ are in bijection with the nested sets of $\Gamma$. These nested sets are our tubings except that the connected components of $\Gamma$ need not be included. Hence there are $2^c$ times as many separable elements as faces of $K(\Gamma)$, where $c$ is the number of connected components. In \cref{sec:nested} we show that the separable elements are exactly the vertices of nested cubes.

\section{Coxeter squares and Coxeter cubes}
\label{sec:cubes}

Coxeter squares were introduced by Dyer \cite[\S 12]{Dyer11II} and Dyer and Wang \cite[\S 3]{Dyer-Wang19}, and Coxeter cubes by the first author \cite[\S 5]{Gimenez25}, following the sketch by Dyer \cite[Example~12.5(1) and \S 14.7]{Dyer11II}. We follow the conventions of \cite{Gimenez25}.

\begin{definition}[{\cite[\S 3]{Dyer-Wang19}, \cite[\S 5]{Gimenez25}}]
\label{def:coxetersquare}
A \emph{Coxeter square} is an ordered quadruple $(w,x,y,z)$ of non-identity elements of $W$ such that
\[
wx = yz \qquad \text{and} \qquad w(\Phi_x) = \Phi_y .
\]
\end{definition}

\begin{figure}[ht]
\centering
\begin{tikzcd}[row sep=large, column sep=large]
\bullet \arrow[d, "w"'] & \bullet \arrow[l, "x"'] \arrow[d, "z"] \\
\bullet & \bullet \arrow[l, "y"]
\end{tikzcd}
\qquad\qquad
\begin{tikzcd}[row sep=large, column sep=large]
w \arrow[d, "w"'] & wx = yz \arrow[l, "x"'] \arrow[d, "z"] \\
e & y \arrow[l, "y"]
\end{tikzcd}
\caption{A Coxeter square $(w,x,y,z)$ as a commutative diagram (left), and with the identity at its terminal vertex (right).}
\label{fig:square}
\end{figure}
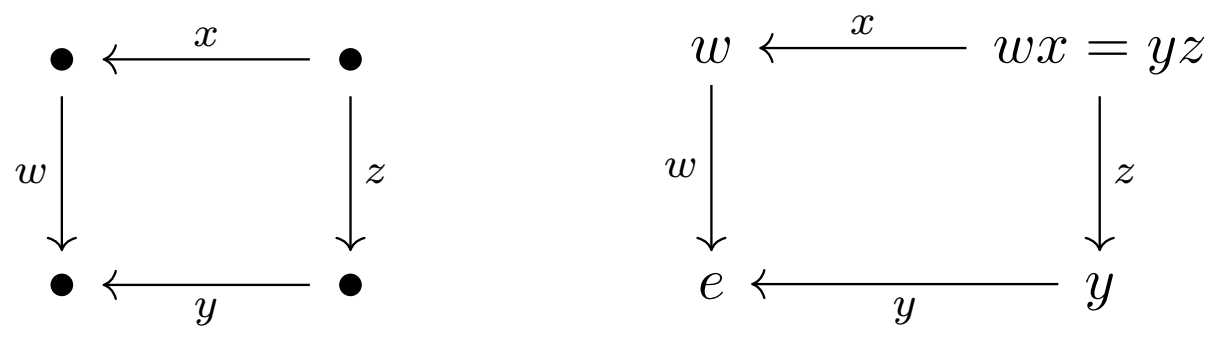

We draw the Coxeter square as the commutative square on the left of \cref{fig:square}. Labels along a directed path are multiplied from right to left, in the order they are applied as functions. A path with labels $g_1, \dots, g_k$ corresponds to the element $g_k \cdots g_1$. The two paths from the initial vertex of the square to its terminal vertex give $wx$ and $yz$. We call $x$ and $z$ the \emph{initial edges} and $w$ and $y$ the \emph{terminal edges} of the square.

If $(w,x,y,z)$ is a Coxeter square then $\ell(wx) = \ell(w) + \ell(x)$ \cite[Proposition~4.1]{Gimenez25}, so \cref{lem:cocycle} gives us the equation
\begin{equation}\label{eqn:squareunion}
\Phi_{wx} = \Phi_w \sqcup w(\Phi_x) = \Phi_w \sqcup \Phi_y .
\end{equation}

A Coxeter square is determined by its terminal edges \cite[Lemma~3.5]{Dyer-Wang19}, but we state it in a slightly different way to fit our needs.

\begin{lemma}
\label{lem:squarecriterion}
Let $w, y \in W$ be non-identity elements. There is a Coxeter square $(w,x,y,z)$ if and only if $\Phi_w \cap \Phi_y = \emptyset$ and $\Phi_w \sqcup \Phi_y$ is an inversion set $\Phi_v$. The square is unique, with $x = w^{-1}v$ and $z = y^{-1}v$.
\end{lemma}

\begin{proof}
If $(w,x,y,z)$ is a Coxeter square, then $\Phi_{wx} = \Phi_w \sqcup \Phi_y$ by \eqref{eqn:squareunion}. Since $v=wx=yz$ is determined by $\Phi_v$, we get $x = w^{-1}v$ and $z = y^{-1}v$.

In the other direction, suppose $\Phi_w$ and $\Phi_y$ are disjoint and their union is an inversion set $\Phi_v$. Then $w \leq_R v$ by \eqref{eqn:weakorder}, so $v = wx$ with $\ell(v) = \ell(w) + \ell(x)$, and \cref{lem:cocycle} gives us $w(\Phi_x) = \Phi_v \setminus \Phi_w = \Phi_y$. Put $z = y^{-1}v$, so that $wx = v = yz$. Neither $x$ nor $z$ is the identity (since $\Phi_w$ and $\Phi_y$ both are nonempty). Hence $(w, x, y, z)$ is a Coxeter square.
\end{proof}

Let $Q_n = \{0,1\}^n$ be the $n$-cube, with each edge oriented so that its varying coordinate increases. We call $\mathbf{0}$ the \emph{initial vertex} and $\mathbf{1}$ the \emph{terminal vertex} of $Q_n$. A two-dimensional face of $Q_n$ has two directed paths from its initial vertex to its terminal vertex. If their labels are $x, w$ and $z, y$ in the order traversed, we say that the face is a Coxeter square when $(w, x, y, z)$ is a Coxeter square. The quadruple $(w,x,y,z)$ is a Coxeter square if and only if $(y,z,w,x)$ is \cite[Lemma~3.4]{Dyer-Wang19} (also clear from \cref{lem:squarecriterion}), so this does not depend on the order of the two paths.

\begin{definition}[{\cite[\S 5]{Gimenez25}}]
\label{def:coxetercube}
A \emph{Coxeter $n$-cube} is a labelling of the edges of $Q_n$ by non-identity elements of $W$ such that every two-dimensional face is a Coxeter square.
\end{definition}

Given a Coxeter $n$-cube, its \emph{terminal edges} are the labels $x_1, \dots, x_n$ of the $n$ edges ending at $\mathbf{1}$. Since every face commutes, the product of the labels along a path between two vertices does not depend on the path. Therefore placing the identity at the terminal vertex $\mathbf{1}$ determines the element at every vertex.

The dimension of a Coxeter cube is at most $r$ \cite[Proposition~5.4]{Gimenez25}, and Coxeter $r$-cubes exist if and only if $W$ is finite \cite[Proposition~5.2 and Corollary~5.11]{Gimenez25}. Moreover we get the following:

\begin{proposition}[{\cite[Propositions~5.7, 5.9 and~5.12]{Gimenez25}}]
\label{prop:cubevertices}
Let $W$ be finite of rank $r$ and let $\C$ be a Coxeter $r$-cube with terminal edges $x_1, \dots, x_r$.
\begin{enumerate}
\item $\C$ is completely determined by its terminal edges.
\item The element at the initial vertex of $\C$ is $w_0$. 
\item The inversion sets $\Phi_{x_1}, \dots, \Phi_{x_r}$ of the terminal edges partition $\Phi^+$.
\end{enumerate}
\end{proposition}

We define a \emph{based cube} to be a Coxeter $r$-cube of a finite Coxeter group $W$ of rank $r$, with the identity at its terminal vertex as above, considered up to permutation of the coordinates of $Q_r$.

\begin{remark}
\label{rem:forgetcoordinates}
By \cref{prop:cubevertices}(1) and \cref{lem:basedcube} below, the coordinates of $Q_r$ carry no information, and we will not use them. For us a based cube is simply a copy of the Boolean lattice inside the right weak order, whose edges are weak order relations (not necessarily the covers).
\end{remark}

\begin{lemma}
\label{lem:basedcube}
Based cubes are exactly the partitions of $\Phi^+$ such that
\begin{itemize}
    \item each block contains a unique simple root, and
    \item the union of any set of blocks is the inversion set of some $w \in W$.
\end{itemize}
\end{lemma}

\begin{proof}
Recall from \cref{subsec:biclosed} that each nonempty inversion set contains at least one simple root. Let $x_1, \dots, x_r$ be the terminal edges of a based cube. By \cref{prop:cubevertices}(3) the inversion sets $\Phi_{x_1}, \dots, \Phi_{x_r}$ partition $\Phi^+$. Each $x_i$ is a non-identity element, so each $\Phi_{x_i}$ contains at least one simple root (exactly one by pigeonhole). The rest of the forward direction follows from \cref{lem:subunion}.

For the backward direction, let $\Phi^+ = B_1 \sqcup \dots \sqcup B_r$ be a partition satisfying the given conditions. For $I \subseteq [r]$ let $u_I$ be the element of $W$ with $\Phi_{u_I} = \bigsqcup_{i \in I} B_i$ (exists by hypothesis). We write $I + i$ for $I \cup \{i\}$. Fix $I$ and distinct $i, j \in [r] \setminus I$. We will show that $u_I, u_{I+i}, u_{I+j}, u_{I+ij}$ are the vertices of a Coxeter square with terminal vertex $u_I$, as in \cref{fig:translatedsquare}. 

Its terminal edges are $w = u_I^{-1} u_{I+i}$ and $y = u_I^{-1} u_{I+j}$, and the product of the labels along either path of the face is $v = u_I^{-1} u_{I+ij}$. Since $\Phi_{u_I} \subseteq \Phi_{u_{I+i}}$, we have $u_I \leq_R u_{I+i}$ by \eqref{eqn:weakorder}, so \cref{lem:cocycle} tells us that $u_I(\Phi_w) = \Phi_{u_{I+i}} \setminus \Phi_{u_I} = B_i$. In the same way $u_I(\Phi_y) = B_j$ and $u_I(\Phi_v) = B_i \sqcup B_j$. Hence $\Phi_w$ and $\Phi_y$ are disjoint with union the inversion set $\Phi_v$, so $(w, w^{-1}v, y, y^{-1}v)$ is a Coxeter square by \cref{lem:squarecriterion}.

\begin{figure}[ht]
\centering
\begin{tikzcd}[row sep=large, column sep=7em]
u_{I+i} \arrow[d, "w = u_I^{-1}u_{I+i}"'] & u_{I+ij} \arrow[l, "w^{-1}v = u_{I+i}^{-1}u_{I+ij}"'] \arrow[d, "y^{-1}v = u_{I+j}^{-1}u_{I+ij}"] \\
u_I & u_{I+j} \arrow[l, "y = u_I^{-1}u_{I+j}"]
\end{tikzcd}
\\[2em]
\begin{tikzcd}[row sep=large, column sep=7em]
\Phi_{u_I} \sqcup B_i \arrow[d, "u_I(\Phi_w) = B_i"'] & \Phi_{u_I} \sqcup B_i \sqcup B_j \arrow[l, "u_{I+i}(\Phi_{w^{-1}v}) = B_j"'] \arrow[d, "u_{I+j}(\Phi_{y^{-1}v}) = B_i"] \\
\Phi_{u_I} & \Phi_{u_I} \sqcup B_j \arrow[l, "u_I(\Phi_y) = B_j"]
\end{tikzcd}
\caption{A square face of the cube constructed in the proof of \cref{lem:basedcube}, with $i, j \notin I$.}
\label{fig:translatedsquare}
\end{figure}
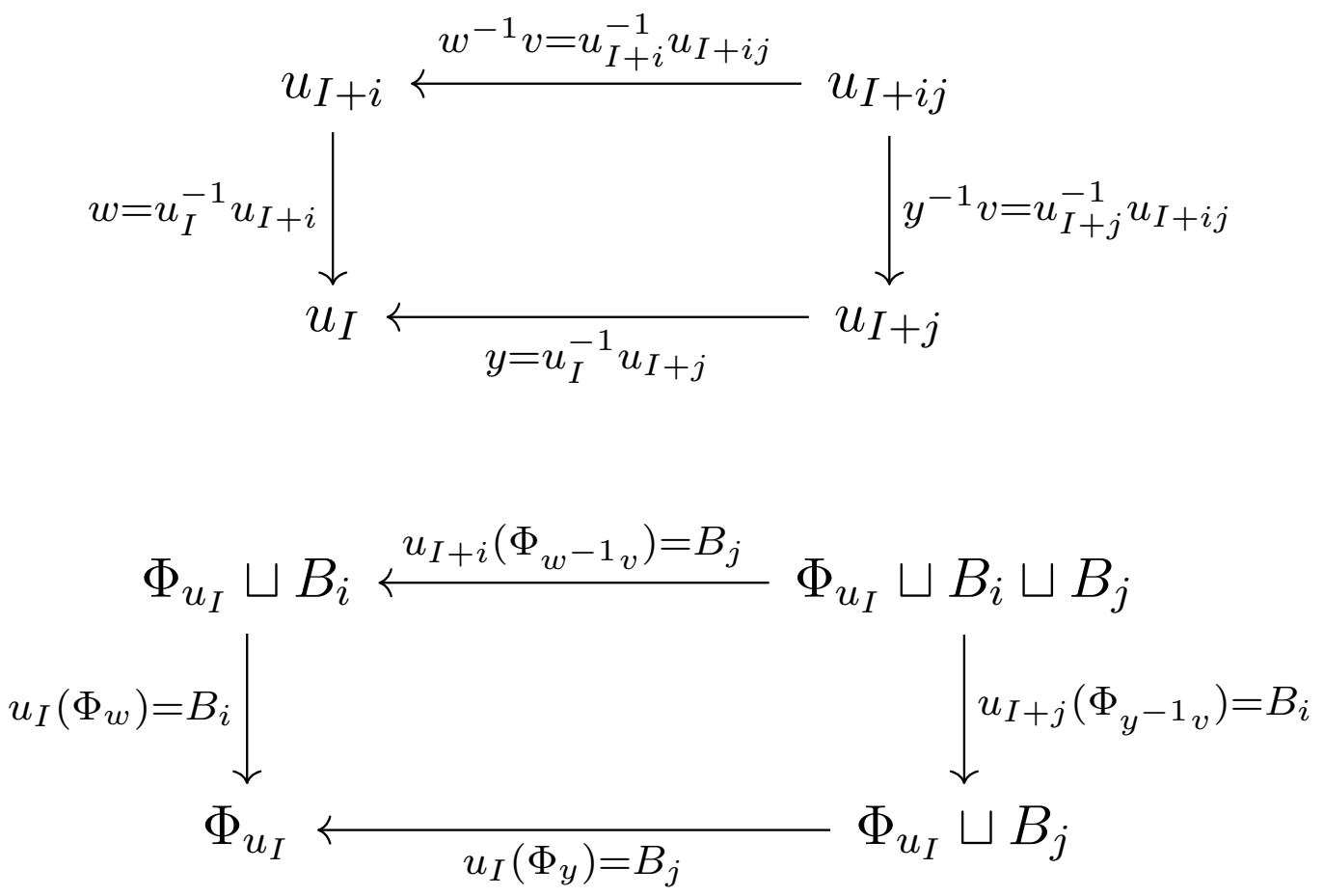

Hence the elements $u_I$ over $I \subseteq [r]$, together with the edges $u_I^{-1}u_{I+i}$ from $u_{I+i}$ to $u_I$ (they are not the identity since $\Phi_{u_I} \subsetneq \Phi_{u_{I+i}}$), form a cube whose square faces are Coxeter squares. Identifying $I \subseteq [r]$ with the $0/1$-vector $e_{[r] \setminus I}$ (a vertex of $Q_r$), this is a Coxeter $r$-cube. Its terminal vertex is $u_\emptyset = e$, so it is a based cube. Its terminal edges are $x_i = u_{\{i\}}$ with $\Phi_{x_i} = B_i$. By \cref{prop:cubevertices}(1) it is the only based cube with these terminal edges.
\end{proof}

If $W$ is reducible, with $S_1, \dots, S_c$ the connected components of $\Gamma$, then $\Phi^+ = \Phi_{S_1}^+ \sqcup \cdots \sqcup \Phi_{S_c}^+$ since supports are connected, and $\Phi_w = \Phi_{w|_{S_1}} \sqcup \cdots \sqcup \Phi_{w|_{S_c}}$ for every $w \in W$ by \eqref{eqn:restriction}. A block $B$ of a based cube is contained in a unique $\Phi_{S_i}^+$: writing $B = \Phi_w$, each nonempty $\Phi_{w|_{S_i}}$ contains a simple root of $S_i$, and $B$ contains only one simple root. It follows from \cref{lem:basedcube} that the based cubes of $W$ are exactly the unions of based cubes of $W_{S_1}, \dots, W_{S_c}$, and we will use this to always reduce $\Gamma$ to the connected case without further comment.

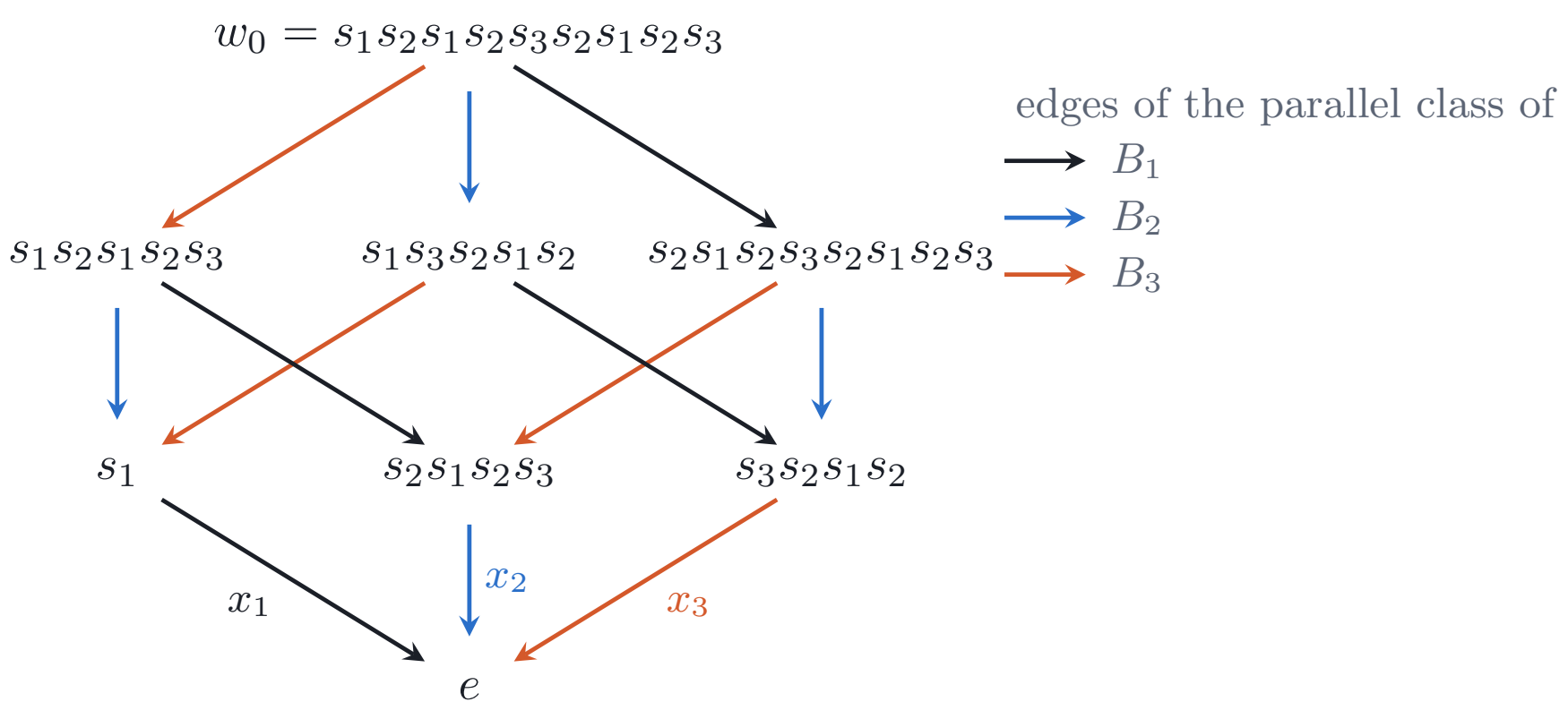
\begin{figure}[ht]
\centering
\input{fig/cube-B3}
\caption{The based cube $\C$ of \cref{ex:B3cube}. The edges of the parallel class corresponding to the block $B_i$ are drawn in the same colour.}
\label{fig:B3cube}
\end{figure}

\begin{example}
\label{ex:B3cube}
Our running example throughout the paper will be $W$ of type $B_3$, with $S = \{s_1, s_2, s_3\}$, $m_{12} = 4$ and $m_{23} = 3$, so that $\Gamma$ is the path $s_1, s_2, s_3$ with the edge $s_1 s_2$ labelled $4$. We write $\alpha_i$ for $\alpha_{s_i}$, and we abbreviate a subset of $S$ by the string of its indices, so that $12$ is the tube $\{s_1, s_2\}$ and $S = 123$. We realize $W$ as the group of signed permutations of $\mathbb{R}^3$ with $\alpha_1 = e_1$, $\alpha_2 = e_2 - e_1$ and $\alpha_3 = e_3 - e_2$. The positive roots are $e_i$ for $1 \le i \le 3$ and $e_j \pm e_i$ for $1 \le i < j \le 3$:
\begin{align*}
\Phi^+ = \{&\alpha_1,\ \alpha_2,\ \alpha_3,\ \alpha_1+\alpha_2,\ \alpha_2+\alpha_3,\ \alpha_1+\alpha_2+\alpha_3, \\
&2\alpha_1+\alpha_2,\ 2\alpha_1+\alpha_2+\alpha_3,\ 2\alpha_1+2\alpha_2+\alpha_3\},
\end{align*}
where $e_i = \alpha_1 + \cdots + \alpha_i$, $e_j - e_i = \alpha_{i+1} + \cdots + \alpha_j$ and $e_j + e_i = 2\alpha_1 + \cdots + 2\alpha_i + \alpha_{i+1} + \cdots + \alpha_j$. Consider the partition $\Phi^+ = B_1 \sqcup B_2 \sqcup B_3$ with
\begin{align*}
B_1 &= \{\alpha_1\} = \{e_1\}, \\
B_2 &= \{\alpha_2,\ \alpha_1+\alpha_2,\ 2\alpha_1+\alpha_2,\ 2\alpha_1+2\alpha_2+\alpha_3\} \\
&= \{e_2 - e_1,\ e_2,\ e_1 + e_2,\ e_2 + e_3\}, \\
B_3 &= \{\alpha_3,\ \alpha_2+\alpha_3,\ \alpha_1+\alpha_2+\alpha_3,\ 2\alpha_1+\alpha_2+\alpha_3\} \\
&= \{e_3 - e_2,\ e_3 - e_1,\ e_3,\ e_1 + e_3\}.
\end{align*}
Each block contains exactly one simple root, and each block is an inversion set. By \cref{prop:chambers}(3), $B_1$, $B_2$ and $B_3$ are the sets of positive roots taking negative values on $\lambda = (-1,2,3)$, $(1,-3,2)$ and $(1,3,-2)$ respectively. Hence every union of blocks is an inversion set by \cref{lem:subunion}, and by \cref{lem:basedcube} the blocks form a based cube $\C$ (drawn in \cref{fig:B3cube}). Its terminal edges are $x_1 = s_1$, $x_2 = s_2s_1s_2s_3$ and $x_3 = s_3s_2s_1s_2$, with $\Phi_{x_i} = B_i$, and its initial vertex is $w_0 = s_1s_2s_1s_2s_3s_2s_1s_2s_3$.

Moving the root $2\alpha_1 + 2\alpha_2 + \alpha_3$ from $B_2$ to $B_3$ gives another partition $\Phi^+ = B_1 \sqcup B_2' \sqcup B_3'$ with
\[
B_2' = \{\alpha_2,\ \alpha_1+\alpha_2,\ 2\alpha_1+\alpha_2\} = \{e_2 - e_1,\ e_2,\ e_1 + e_2\}, \qquad B_3' = B_3 \sqcup \{e_2 + e_3\},
\]
which is again a based cube $\C'$. We will see in \cref{ex:B3nested} that $\C'$ is a nested cube while $\C$ is not. In \cref{ex:B3all} we will see where both cubes sit on the graph associahedron (one at a vertex and one along an edge).
\end{example}

\section{Nested cubes and the graph associahedron}
\label{sec:nested}

In this section we study \emph{nested cubes}. These are the based cubes where all blocks are parabolic. We show that they are in bijection with the vertices of the graph associahedron $K(\Gamma)$ of the Coxeter diagram, and that the elements of $W$ appearing as vertices of nested cubes are exactly the separable elements of Gaetz and Gao.

\begin{definition}
\label{def:parabolic}
A subset $B$ of $\Phi^+$ is \emph{parabolic} if
\[
B = B(T,p) := \{\beta \in \Phi^+ : p \in \supp(\beta) \subseteq T\}
\]
for a tube $T$ of $\Gamma$ and a node $p \in T$.
\end{definition}

We call a based cube a \emph{nested cube} if all of its blocks are parabolic. 

\begin{lemma}
\label{lem:parabolicinversion}
Every parabolic block $B(T,p)$ is an inversion set.
\end{lemma}
\begin{proof}
Let $T_1, \dots, T_k$ be the connected components of $T \setminus \{p\}$. A root in $\Phi_T^+$ either contains $p$ in its support, or has its support in some $T_i$. In the first case we get $B(T,p)$, in the second case we get $\Phi^+_{T \setminus \{p\}} = \Phi_u$ with $u = w_0(T \setminus \{p\})$. Hence $B(T,p)$ is the complement of $\Phi_u$ in $\Phi_T^+$, which is an inversion set by \eqref{eqn:paraboliccomplement}.
\end{proof}

\begin{theorem}
\label{thm:nestedcubesgamma}
There is a bijection between nested cubes and vertices of $K(\Gamma)$.
\end{theorem}
\begin{proof}
We first show that a maximal tubing $N$ gives a nested cube. Since $N$ is maximal, every $T \in N$ has a pivot $p(T)$. For each $T \in N$ take the block $B(T,p(T))$, which is an inversion set by \cref{lem:parabolicinversion}. These blocks are pairwise disjoint, since any two tubes of $N$ are either disjoint (in this case the supports of the roots in the two blocks are disjoint), or nested $T' \subsetneq T$ (in this case we have $p(T) \notin T'$). We next show that the blocks cover $\Phi^+$. Take any $\beta \in \Phi^+$ and let $T$ be the smallest tube of $N$ containing $\supp(\beta)$, which exists since $S \in N$. If $p(T) \notin \supp(\beta)$, then $\supp(\beta)$ is a connected subset of $T \setminus \{p(T)\}$. Hence $\supp(\beta)$ lies in a child $T'$ of $T$ as well, contradicting the choice of $T$. Therefore we get that $\beta \in B(T,p(T))$.
Since $|N| = r$, \cref{lem:subunion} and \cref{lem:basedcube} give us that these blocks form a based cube.

For the other direction, let a nested cube consist of parabolic blocks $B(T_i,p_i)$, $1 \le i \le r$. We first show that the $T_i$ are pairwise compatible. If $T_i$ and $T_j$ are disjoint and adjacent, or overlap without being nested, then $T_i \cup T_j$ is a tube. So there is a root with support $T_i \cup T_j$, and it lies in some block $B(T,p)$ of the cube, with $p \in T_i \cup T_j \subseteq T$. Without loss of generality $p \in T_i$. A root with support $T_i$ then lies in both $B(T_i,p_i)$ and $B(T,p)$, so these blocks coincide and $T = T_i$. This contradicts the assumption that $T_i$ and $T_j$ are not nested.

The same reasoning shows that the $T_i$'s are distinct, since if $T_i = T_j$ with $i \neq j$ a root with support $T_i$ would lie in both blocks. Hence the $T_i$'s are $r$ distinct pairwise compatible tubes, so they form a maximal tubing $N$. Finally $p_i = p(T_i)$, since if $p_i$ is in a child $T_j$ of $T_i$ in $N$, a root with support $T_j$ would lie in both $B(T_j,p_j)$ and $B(T_i,p_i)$, giving us a contradiction. Hence the two constructions we have described are inverse maps to each other.
\end{proof}

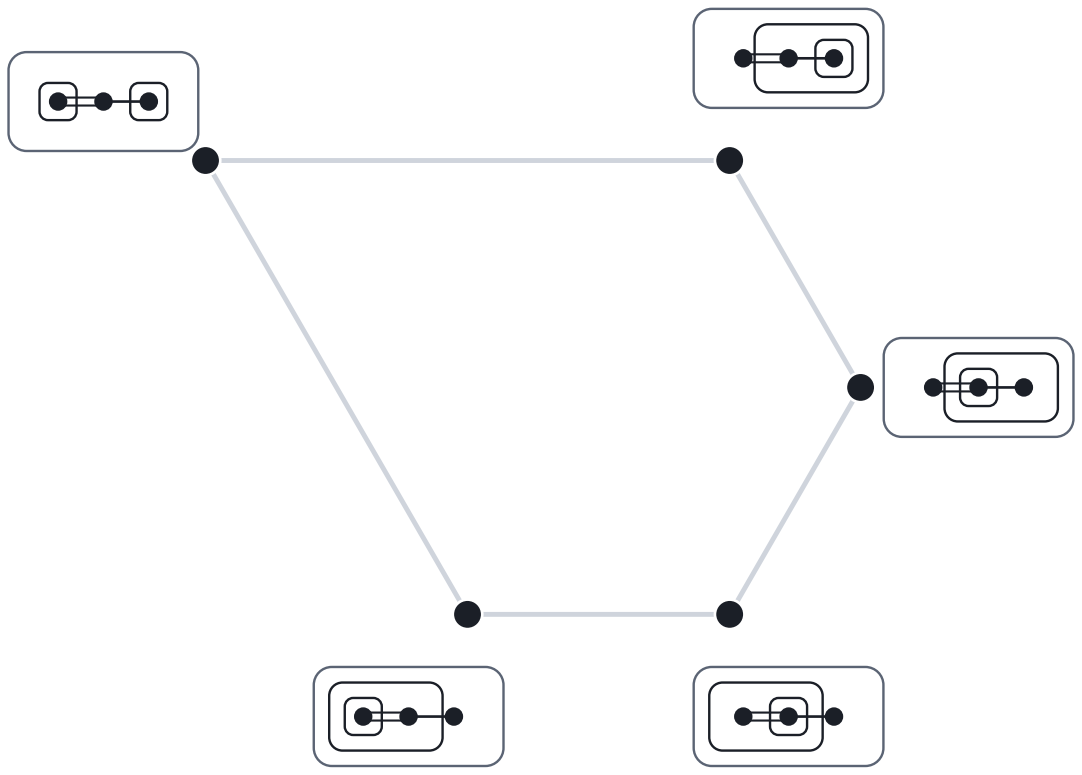
\begin{figure}[ht]
\centering
\input{fig/pentagon-B3-nested}
\caption{The graph associahedron $K(\Gamma)$ for $\Gamma$ of type $B_3$ (a pentagon), with the maximal tubing at each vertex. In the Coxeter diagrams the nodes are $s_1, s_2, s_3$ from left to right.}
\label{fig:B3pentagon}
\end{figure}

\begin{example}
\label{ex:B3nested}
Let $W$ be our $B_3$ running example from \cref{ex:B3cube}. Since $\Gamma$ is a path on three nodes, $K(\Gamma)$ is a pentagon. Its five vertices are the maximal tubings
\[
\{S, 12, 1\}, \quad \{S, 12, 2\}, \quad \{S, 23, 2\}, \quad \{S, 23, 3\}, \quad \{S, 1, 3\},
\]
shown in \cref{fig:B3pentagon}. The parabolic blocks of $\Gamma$ are
\begin{align*}
B(1,1) &= \{\alpha_1\}, & B(12,1) &= \{\alpha_1,\ \alpha_1+\alpha_2,\ 2\alpha_1+\alpha_2\}, \\
B(2,2) &= \{\alpha_2\}, & B(12,2) &= \{\alpha_2,\ \alpha_1+\alpha_2,\ 2\alpha_1+\alpha_2\}, \\
B(3,3) &= \{\alpha_3\}, & B(23,2) &= \{\alpha_2,\ \alpha_2+\alpha_3\}, \\
& & B(23,3) &= \{\alpha_3,\ \alpha_2+\alpha_3\},
\end{align*}
together with the blocks $B(S,p) = \Phi^+ \setminus \Phi^+_{S \setminus \{p\}}$ of full support,
\begin{align*}
B(S,1) &= \Phi^+ \setminus \{\alpha_2,\ \alpha_3,\ \alpha_2+\alpha_3\}, \\
B(S,2) &= \Phi^+ \setminus \{\alpha_1,\ \alpha_3\}, \\
B(S,3) &= \Phi^+ \setminus \{\alpha_1,\ \alpha_2,\ \alpha_1+\alpha_2,\ 2\alpha_1+\alpha_2\},
\end{align*}
with six, seven and five elements respectively. The five nested cubes are
\begin{align*}
\{S, 12, 1\} &\longmapsto B(1,1) \sqcup B(12,2) \sqcup B(S,3), \\
\{S, 12, 2\} &\longmapsto B(2,2) \sqcup B(12,1) \sqcup B(S,3), \\
\{S, 23, 2\} &\longmapsto B(2,2) \sqcup B(23,3) \sqcup B(S,1), \\
\{S, 23, 3\} &\longmapsto B(3,3) \sqcup B(23,2) \sqcup B(S,1), \\
\{S, 1, 3\} &\longmapsto B(1,1) \sqcup B(3,3) \sqcup B(S,2).
\end{align*}
The first one is the cube $\C'$ from \cref{ex:B3cube}, since $B_2' = B(12,2)$ and $B_3' = B(S,3)$. The cube $\C$ of \cref{ex:B3cube} is not nested: a parabolic block containing $\alpha_2$ and the root $2\alpha_1 + 2\alpha_2 + \alpha_3$ of full support would have to be $B(S,2)$, which has seven elements, but $B_2$ has only four.
\end{example}

For finite Weyl groups the equivalence of (1) and (3) below essentially follows from the bijection of \cite[Theorem~7]{Gaetz-Gao20split} between separable
elements and nested sets of $\Gamma$. We state and prove it for all finite Coxeter
groups.

\begin{theorem}\label{thm:nestedvertices}
For $w \in W$ the following are equivalent.
\begin{enumerate}
    \item $w$ is a vertex of some nested cube.
    \item $\Phi_w$ and $\Phi^+ \setminus \Phi_w$ both can be partitioned into parabolic blocks.
    \item $w$ is separable.
\end{enumerate}
\end{theorem}
\begin{proof}
We first show the equivalence of (1) and (2). If $w$ is a vertex of a nested cube, then $\Phi_w$ is a union of blocks of the cube by \cref{lem:basedcube}, and these blocks are parabolic. They partition $\Phi_w$ and the remaining blocks partition $\Phi^+ \setminus \Phi_w$.

In the opposite direction, suppose $\Phi_w$ and $\Phi^+ \setminus \Phi_w$ are partitioned into parabolic blocks. A parabolic block $B(T,p)$ contains exactly one simple root $\alpha_p$, and every simple root lies in one of the blocks. Therefore, there are exactly $r$ blocks $B_1, \dots, B_r$, with $\Phi_w = B_1 \sqcup \cdots \sqcup B_k$. These are inversion sets by \cref{lem:parabolicinversion}, so every union of blocks is an inversion set by \cref{lem:subunion}. By \cref{lem:basedcube} the blocks form a based cube, which is nested since all its blocks are parabolic. Then $w$ is its vertex corresponding to $B_1 \sqcup \cdots \sqcup B_k$.

We now show the equivalence of (2) and (3), by induction on the rank $r$. In the base case $r = 1$, both elements of $W$ are separable by definition, and $\Phi^+ = \{\alpha_s\}$ is a single parabolic block, so (2) holds for both as well. Now assume that the claim holds in smaller rank. If $W$ is reducible, parabolic blocks stay within a single component, so both (2) and (3) reduce to the components. Hence we only need to consider the case when $W$ is irreducible.

Recall $B(S,p) = \{\beta \in \Phi^+ : p \in \supp(\beta)\}$, so $p$ is a full or empty pivot of $w$ if and only if $B(S,p)$ lies in $\Phi_w$ or in $\Phi^+ \setminus \Phi_w$. Starting from (2), a root of full support lies in some block $B(S,p)$ where $p$ is a full or empty pivot. The remaining blocks partition $\Phi_w \cap \Phi^+_{S \setminus \{p\}}$ and its complement in $\Phi^+_{S \setminus \{p\}}$ into parabolic blocks of $\Gamma \setminus \{p\}$, so $w|_{S \setminus \{p\}}$ is separable by the induction hypothesis. Starting from (3), given a full or empty pivot $p$, the induction hypothesis partitions $\Phi_{w|_{S \setminus \{p\}}}$ and its complement in $\Phi^+_{S \setminus \{p\}}$ into parabolic blocks, which are parabolic blocks of $\Gamma$ as well. Adding $B(S,p)$ to $\Phi_{w|_{S \setminus \{p\}}}$ (when $p$ is a full pivot) or to its complement in $\Phi^+_{S \setminus \{p\}}$ (when $p$ is an empty pivot) gives (2).
\end{proof}

\section{Non-nested cubes on cubical faces}
\label{sec:nonnested}

In the previous section we identified the nested cubes with the vertices of the graph associahedron. In this section we describe the remaining cubes. We show that every based cube lies in the interior of a cubical face of $K(\Gamma)$, and that the number of cubes on an edge comes from a rank two interval of the intersection lattice of the Coxeter arrangement.

The following structure theorem is due to Dimitrov, Gigliotti, Ossip, Paquette and Wehlau \cite{Dimitrov-etal25} in the case of finite Weyl groups. They study decompositions of $\Phi^+$ into inversion sets for a finite crystallographic root system and call the decompositions into $r$ inversion sets \emph{fine decompositions}. By \cref{lem:basedcube} these are exactly our based cubes.

Proposition 10.1 of \cite{Dimitrov-etal25} shows that in every based cube there is $S' \subseteq S$ with $|S'| = r-1$ or $r-2$ such that all blocks but one or two are contained in $\Phi_{S'}$. Corollary 9.2 of \cite{Dimitrov-etal25} shows that these remaining blocks are exactly the blocks of full support.

\begin{theorem}
\label{thm:structure}
Let $W$ be an irreducible finite Coxeter group. Every based cube is exactly one of the following two cases:
\begin{itemize}
    \item One block is $\Phi^+ \setminus \Phi^+_{S \setminus \{p\}}$ for a node $p$, or
    \item Two blocks partition $\Phi^+ \setminus \Phi^+_{S \setminus \{p,q\}}$ for nodes $p \neq q$, and both blocks have full support.
\end{itemize}
In both cases the remaining blocks form a based cube of $W_{S'}$, where $S' = S \setminus \{p\}$ or $S' = S \setminus \{p,q\}$ respectively.
\end{theorem}

\begin{proof}
As explained above, for finite Weyl groups this is \cite[Corollary 9.2 and Proposition 10.1]{Dimitrov-etal25}.

For $H_3$ and $H_4$ the statement was verified by enumerating all based cubes ($16$ and $118$ respectively) and checking the recursion on each of them \cite{GO26code}.

For $I_2(m)$ with $S = \{p,q\}$, take $J = \emptyset$, so that $X_J = V$ and $U = \Phi^+$. Write the positive roots as $\beta_1 = \alpha_p, \dots, \beta_m = \alpha_q$ in increasing order of the ratio $\rho$ (distinct roots have distinct ratios). By \cref{lem:plane} and the discussion before it, the inversion sets are exactly the sets $B_\lambda$ for $\lambda$ on no hyperplane of $\A$, and these are $\emptyset$, $\Phi^+$ and the initial and final segments of $\{\beta_1,\dots,\beta_m\}$. Hence the based cubes are exactly the pairs
$$\{\beta_1, \dots, \beta_k\}, \{\beta_{k+1}, \dots, \beta_m\}, \qquad 1 \le k \le m-1.$$
For $k=1$ or $k=m-1$ we are in the first case. For other values of $k$ both blocks have full support and we are in the second case.

Regardless of the type, the two cases cannot occur at the same time: the block $B(S,p) = \Phi^+ \setminus \Phi^+_{S \setminus \{p\}}$ of the first case contains every root of full support (support $S$), so no two disjoint blocks can both have full support.
\end{proof}

In type $A$ only the first case occurs, and the statement was proved by the first author \cite[Theorem~5.20]{Gimenez25} using Brink--Howlett generators \cite{Brink-Howlett99}. 

Let $J = S \setminus \{p,q\}$, let $L_1, \dots, L_\sigma$ be as in \cref{subsec:ranktwo}, and put $U = \Phi^+ \setminus \Phi_J^+$.

\begin{lemma}\label{lem:prefix}
Let $B \subseteq U$ with $\alpha_p \in B$ and $\alpha_q \notin B$. Then $B$ and $U \setminus B$ are both inversion sets if and only if $B = L_1 \sqcup \cdots \sqcup L_k$ for some $1 \leq k \leq \sigma - 1$.
\end{lemma}

\begin{proof}
Suppose $B = L_1 \sqcup \cdots \sqcup L_k$ with $1 \le k \le \sigma - 1$. Let $Q := X_J \cap H_{\alpha_p}^- \cap H_{\alpha_q}^+$, an open quadrant of $X_J$ where the boundary rays are on $H_{\alpha_p}$ and $H_{\alpha_q}$. For $\beta \in U$ with $0 < \rho(\beta) < \infty$ the line $H_\beta \cap X_J$ meets $Q$ in a ray, and these rays are ordered by $\rho$ from the boundary ray on $H_{\alpha_p}$ ($\rho = 0$) to the one on $H_{\alpha_q}$ ($\rho = \infty$). Pick $\lambda \in Q$ strictly between the rays of $L_k$ and $L_{k+1}$ ($\rho_k < -\lambda_p/\lambda_q < \rho_{k+1}$). Then $B_\lambda = B$ and $\lambda$ lies on no $H_\beta$ with $\beta \in U$, so \cref{lem:plane}(1) gives us that $B$ and $B \sqcup \Phi_J^+$ are inversion sets. Moreover from \eqref{eqn:complement} we get that $U \setminus B = \Phi^+ \setminus (B \sqcup \Phi_J^+)$ is an inversion set as well.

For the other direction, suppose that $B$ and $U \setminus B$ are both inversion sets. Again from \eqref{eqn:complement} we get that $B \sqcup \Phi_J^+ = \Phi^+ \setminus (U \setminus B)$ is an inversion set as well. From \cref{lem:plane}(2) we get $B = B_\lambda$ for some $\lambda \in X_J$. Since $\alpha_p \in B_\lambda$ and $\alpha_q \notin B_\lambda$ we have $\lambda_p < 0 < \lambda_q$, so $B_\lambda$ partitions as $L_1 \sqcup \cdots \sqcup L_k$.
\end{proof}

\begin{figure}[ht]
\centering
\resizebox{0.46\textwidth}{!}{\input{fig/flat-B3}}\hfill
\resizebox{0.46\textwidth}{!}{\input{fig/hemisphere-B3}}
\caption{Left figure shows the lines $H_\beta \cap X_J$ for $\beta \in U$ of \cref{ex:B3flat} in the plane $X_J = \{x_1 = 0\}$ with coordinates $(x_2,x_3)$. The cones corresponding to $B_2$ and $B_3$ are shaded. Right figure is the hemisphere $x_1 > 0$ of the unit sphere in $V$ viewed along $e_1$, with $X_J \cap S^2$ as its rim.}
\label{fig:B3flat}
\end{figure}
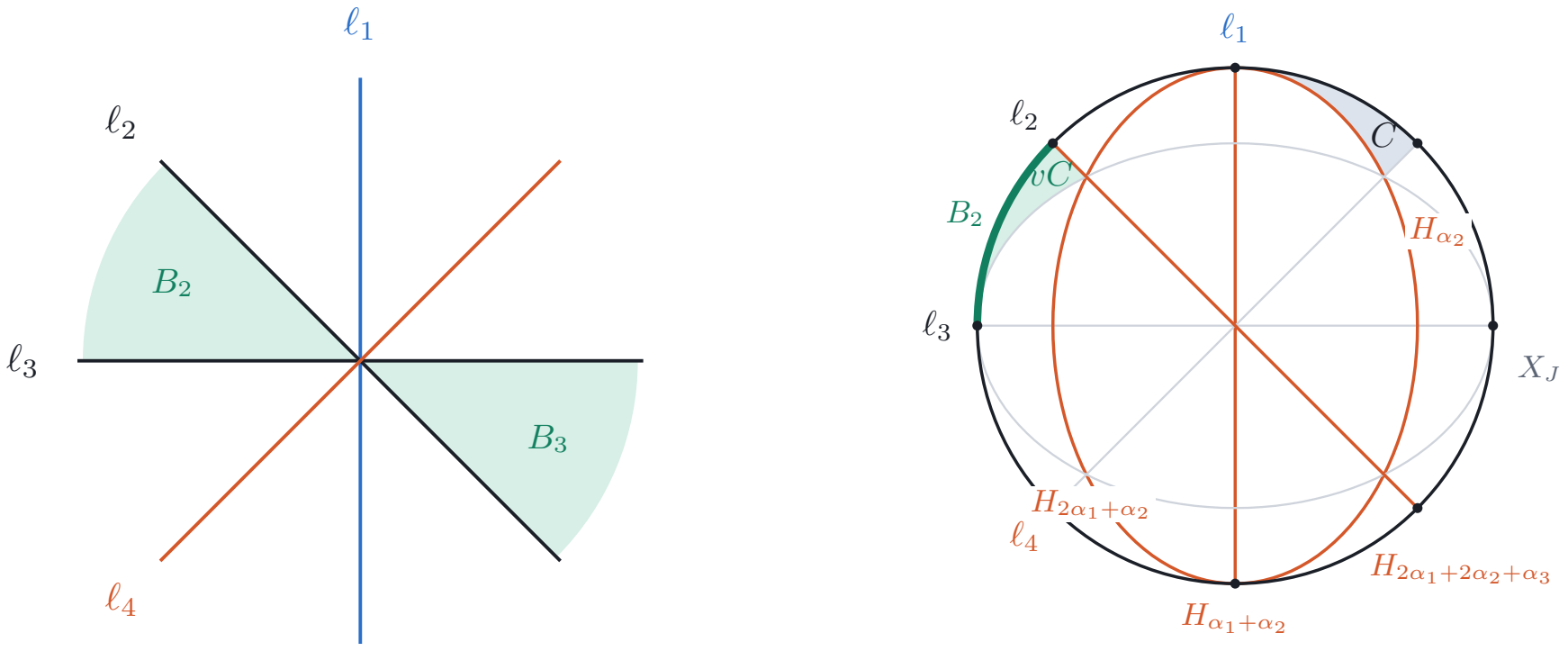

\begin{example}
\label{ex:B3flat}
Let $W$ be our $B_3$ running example from \cref{ex:B3cube}, and take $J = \{s_1\}$, $p = s_2$ and $q = s_3$. Then $X_J = H_{\alpha_1} = \{x_1 = 0\}$ (we identify this with the $(x_2, x_3)$-plane), $\Phi_J^+ = \{\alpha_1\}$, and $U = \Phi^+ \setminus \Phi_J^+$ consists of the remaining eight roots. For $\lambda \in X_J$ and $\beta = c_1\alpha_1 + c_2\alpha_2 + c_3\alpha_3$ we have $(\beta, \lambda) = c_2(\alpha_2,\lambda) + c_3(\alpha_3,\lambda)$, so the line $H_\beta \cap X_J$ only depends on the ratio $\rho(\beta) = c_3/c_2$. The ratios occurring in $U$ are $0$, $1/2$, $1$ and $\infty$, so we get $\sigma = 4$. The classes are
\begin{align*}
L_1 &= \{\alpha_2,\ \alpha_1+\alpha_2,\ 2\alpha_1+\alpha_2\} = \{e_2 - e_1,\ e_2,\ e_1 + e_2\}, \\
L_2 &= \{2\alpha_1+2\alpha_2+\alpha_3\} = \{e_2 + e_3\}, \\
L_3 &= \{\alpha_2+\alpha_3,\ \alpha_1+\alpha_2+\alpha_3,\ 2\alpha_1+\alpha_2+\alpha_3\} = \{e_3 - e_1,\ e_3,\ e_1 + e_3\}, \\
L_4 &= \{\alpha_3\} = \{e_3 - e_2\},
\end{align*}
and the lines $\ell_i := H_\beta \cap X_J$ for $\beta \in L_i$ are $\ell_1 = \{x_2 = 0\}$, $\ell_2 = \{x_2 + x_3 = 0\}$, $\ell_3 = \{x_3 = 0\}$ and $\ell_4 = \{x_3 = x_2\}$, shown on the left of \cref{fig:B3flat}. They cut $X_J$ into eight open cones and $B_\lambda$ is constant on each of them. By \cref{lem:prefix}, the sets $B \subseteq U$ with $\alpha_2 \in B$ and $\alpha_3 \notin B$ such that $B$ and $U \setminus B$ are inversion sets, are the following:
\[
L_1 = B(12, 2), \qquad L_1 \sqcup L_2 = B_2, \qquad L_1 \sqcup L_2 \sqcup L_3 = B(S, 2),
\]
where $B_2$ is the block of the cube $\C$ in \cref{ex:B3cube}. In particular $\lambda = (-3,2) \in X_J$ gives $B_\lambda = L_1 \sqcup L_2 = B_2$, and the two elements $v$ with $\Phi_v \setminus \{\alpha_1\} = B_2$ are $s_2s_1s_2s_3$ and $s_1s_2s_1s_2s_3$, with inversion sets $B_2$ and $B_1 \sqcup B_2$ respectively. The right side of \cref{fig:B3flat} shows the chamber $vC$ for $v = s_2s_1s_2s_3$ as the spherical triangle attached to the arc of the rim where we have $B_\lambda = B_2$, and the four red circles $H_\beta \cap S^2$ with $\beta \in B_2$ separate $vC$ from $C$. The chamber for $s_1s_2s_1s_2s_3$ is its mirror image in $X_J$ on the back hemisphere. The opposite point $-\lambda$ gives $B_{-\lambda} = L_3 \sqcup L_4 = B_3$ and the elements $s_3s_2s_1s_2$ and $s_1s_3s_2s_1s_2$ with inversion sets $B_3$ and $B_1 \sqcup B_3$ respectively. These four elements are vertices of $\C$ in \cref{fig:B3cube}.
\end{example}

\subsection{The main result}

Let $T$ be a tube and let $p, q \in T$ be distinct nodes. Applying \cref{subsec:ranktwo} to the group $W_T$ with its root system $\Phi_T$ and $J = T \setminus \{p,q\}$, we get the classes $L_1(T,p,q), \dots, L_\sigma(T,p,q)$ of $\Phi_T^+ \setminus \Phi_J^+$, where two roots are in the same class when they have the same ratio $c_q(\beta)/c_p(\beta)$, and the classes are ordered by increasing ratio. We use $\sigma(T,p,q)$ to denote their count. This count also equals the number of atoms of the rank two interval $[X_J, X_T]$ in $\LL(\A)$. The atoms are the flats $H_\beta \cap X_J$ with $\beta \in \Phi_T^+ \setminus \Phi_J^+$, and $H_\beta \cap X_J = H_\gamma \cap X_J$ if and only if $\beta$ and $\gamma$ have the same ratio. We have $\sigma(T,p,q) \geq 3$, since $\alpha_p$, $\alpha_q$ and any root with support $T$ have ratios $0$, $\infty$ and a finite positive value respectively. In type $A$ we always have $\sigma(T,p,q) = 3$.

Recall from \cref{subsec:assoc} that a tubing is cubical if every shell has size at most two. Given a based cube $\C$ we obtain a cubical tubing $N(\C)$ in a recursive way. If $\Gamma$ is disconnected, $N(\C)$ is the union of the cubical tubings coming from the based cubes of the components. If $\Gamma$ is connected, then from \cref{thm:structure} either one block of $\C$ is the parabolic block $B(S,p) = \Phi^+ \setminus \Phi^+_{S'=S \setminus \{p\}}$, or two blocks partition $\Phi^+ \setminus \Phi^+_{S'=S \setminus \{p,q\}}$. In both cases the remaining blocks form a based cube of $W_{S'}$. By the discussion following \cref{lem:basedcube}, this is a product of based cubes of the connected components of $S'$, and we recurse on each of them. The resulting tubes form a tubing where the children of $S$ are the components of $S'$, such that the shell of $S$ is $\{p\}$ or $\{p,q\}$. From the recursive construction $N(\C)$ is cubical. Notice that $N(\C)$ is maximal if and only if the first case occurs at every step through the recursion. Hence a cube is nested if and only if $N(\C)$ is maximal.

For a cubical tubing $N$ let $N_2 \subseteq N$ be the set of tubes with a two-element shell. For each $T \in N_2$ fix an order $(p_T, q_T)$ of its shell, so that $\mathcal{F}_N$ is a cube of dimension $|N_2|$ with one edge for each $T \in N_2$.

Now let $\C$ be a based cube with $N(\C) = N$ and let $T \in N_2$, with $J = T \setminus \{p_T, q_T\}$. In the construction of $N(\C)$, the second case of \cref{thm:structure} occurred at $T$: two blocks $B \ni \alpha_{p_T}$ and $B' \ni \alpha_{q_T}$ of $\C$ partition $\Phi_T^+ \setminus \Phi_J^+$, and the remaining blocks of $\C$ inside $\Phi_T^+$ form a based cube of $W_J$. Since $B$ and $B'$ are inversion sets, by \cref{lem:prefix} applied to $W_T$ there is a unique integer $k_T(\C)$ with
\[
B = L_1(T,p_T,q_T) \sqcup \cdots \sqcup L_{k_T(\C)}(T,p_T,q_T), \qquad 1 \le k_T(\C) \le \sigma(T,p_T,q_T) - 1 .
\]
We define $\pi(\C)$ to be the point of the cubical face $\mathcal{F}_N \cong [0,1]^{N_2}$ with coordinates
\[
\pi(\C)_T = \frac{k_T(\C) - 1}{\sigma(T,p_T,q_T) - 2}, \qquad T \in N_2 ,
\]
so that the coordinate of $T$ is $0$ or $1$ at the endpoints (with $0$ at the endpoint $N \cup \{T_{p_T}\}$, as in \cref{subsec:assoc}). This defines a map $\pi$ from the set of based cubes of $W$ to $K(\Gamma)$.

\begin{theorem}
\label{thm:cubicaltubing}
The map $\pi$ is injective, and it has the following properties.
\begin{enumerate}
\item On nested cubes, $\pi$ is the bijection of \cref{thm:nestedcubesgamma} onto the vertices of $K(\Gamma)$. Every other based cube $\C$ is sent to the interior of the cubical face $\mathcal{F}_{N(\C)}$.
\item An edge of $K(\Gamma)$ with two-element shell $\{p,q\}$ in the tube $T$ contains exactly $\sigma(T,p,q) - 3$ equally spaced based cubes in its interior.
\item A cubical face $\mathcal{F}_N$ is the product of its edges, and the set of points $\pi(\C)$ on it is the product of the sets of points on the edges.
\end{enumerate}
In particular, the number of based cubes in the interior of a cubical face $\mathcal{F}_N$ is
\[
\prod_{T \in N_2} \bigl(\sigma(T,p_T,q_T) - 3\bigr),
\]
and every based cube lies in the interior of exactly one cubical face.
\end{theorem}

\begin{proof}
We show that for every cubical tubing $N$ the map $\C \mapsto (k_T(\C))_{T \in N_2}$ is a bijection from the set of based cubes with $N(\C) = N$ onto the set of tuples of integers $(k_T)_{T \in N_2}$ with $2 \le k_T \le \sigma(T,p_T,q_T) - 2$. 

We proceed by induction on $|S|$. In the base case $|S| = 1$ the statement is trivial. Assume the claim holds for smaller $|S|$. If $\Gamma$ is disconnected, a based cube is a product of based cubes of the components and a cubical tubing is a union of cubical tubings of the components, so the claim follows. Hence we may assume that $\Gamma$ is connected.

Let $\C$ be a based cube with $N(\C) = N$ and apply \cref{thm:structure} at $S$ with $S' = S \setminus \{p\}$ or $S' = S \setminus \{p,q\}$. In both cases the remaining blocks form a based cube of $W_{S'}$, and any based cube of $W_{S'}$ together with $|S \setminus S'|$ nonempty inversion sets partitioning $\Phi^+ \setminus \Phi_{S'}^+$ forms a based cube of $W$ thanks to \cref{lem:subunion,lem:basedcube}.

In the first case $S$ has shell $\{p\}$ and $N_2$ is unchanged, so the claim follows from the induction hypothesis. In the second case $S$ has shell $\{p,q\}$ and drops from $N_2$, so by the induction hypothesis the based cubes of $W_{S'}$ with tubing $N \setminus \{S\}$ are in bijection with the tuples $(k_T)_{T \in N_2 \setminus \{S\}}$. All that remains for us is to identify the pairs $B \ni \alpha_p$, $B' \ni \alpha_q$ partitioning $U = \Phi^+ \setminus \Phi_{S'}^+$ with the integers $2 \le k_S \le \sigma(S,p,q) - 2$.

By \cref{lem:prefix} we have $B = L_1 \sqcup \cdots \sqcup L_k$ and $B' = L_{k+1} \sqcup \cdots \sqcup L_\sigma$ for some $1 \le k \le \sigma - 1$. Since $L_1$ is the set of roots whose support contains $p$ but not $q$, we have $B(S,q) = U \setminus L_1$, and similarly $B(S,p) = U \setminus L_\sigma$. Hence for $k = 1$ or $k = \sigma - 1$ we are in the first case of \cref{thm:structure} and $N(\C) \neq N$. For $2 \le k \le \sigma - 2$ neither $B$ nor $B'$ is of the form $B(S,\cdot)$, so we are in the second case with tubing $N$. This completes the induction.

Finally, the nested cubes are those with $N_2 = \emptyset$, so they are placed at the vertices of $K(\Gamma)$ (consistent with \cref{thm:nestedcubesgamma}). Every other based cube $\C$ has all coordinates of $N_2$ strictly between $0$ and $1$, so it lies in the interior of $\mathcal{F}_{N(\C)}$. Distinct cubes with the same tubing are assigned distinct points, giving us injectivity and the first item. The third item is the bijection with tuples just shown (this gives the bijection on the interior, and applying it to the faces of $\mathcal{F}_N$ gives the bijection on the whole face). The second item is the special case when $|N_2| = 1$.
\end{proof}

The proof describes all the based cubes on a cubical face $\mathcal{F}_N$. Along the edge of $T \in N_2$ the block containing $\alpha_{p_T}$ is $L_1 \sqcup \cdots \sqcup L_{k_T}$, so consecutive based cubes differ by moving a single class $L_i$ from the block containing $\alpha_{q_T}$ to the block containing $\alpha_{p_T}$.

\begin{corollary}
\label{cor:count}
The number of based cubes of $W$ is
\[
\sum_{N} \prod_{T \in N_2} \bigl(\sigma(T,p_T,q_T) - 3\bigr),
\]
the sum over all cubical tubings $N$ of $\Gamma$. In type $A$ only the maximal tubings contribute, and the based cubes are exactly the nested cubes.
\end{corollary}

\begin{figure}[ht]
\centering
\input{fig/pentagon-B3-twisted}
\caption{The nine based cubes of $B_3$ as points of the pentagon $K(\Gamma)$, as in \cref{ex:B3all}. The vertices are the five nested cubes of \cref{fig:B3pentagon}, each edge is labelled by its value of $\sigma$. The four edges with $\sigma = 4$ carry one non-nested cube each, at their midpoints.}
\label{fig:B3twisted}
\end{figure}
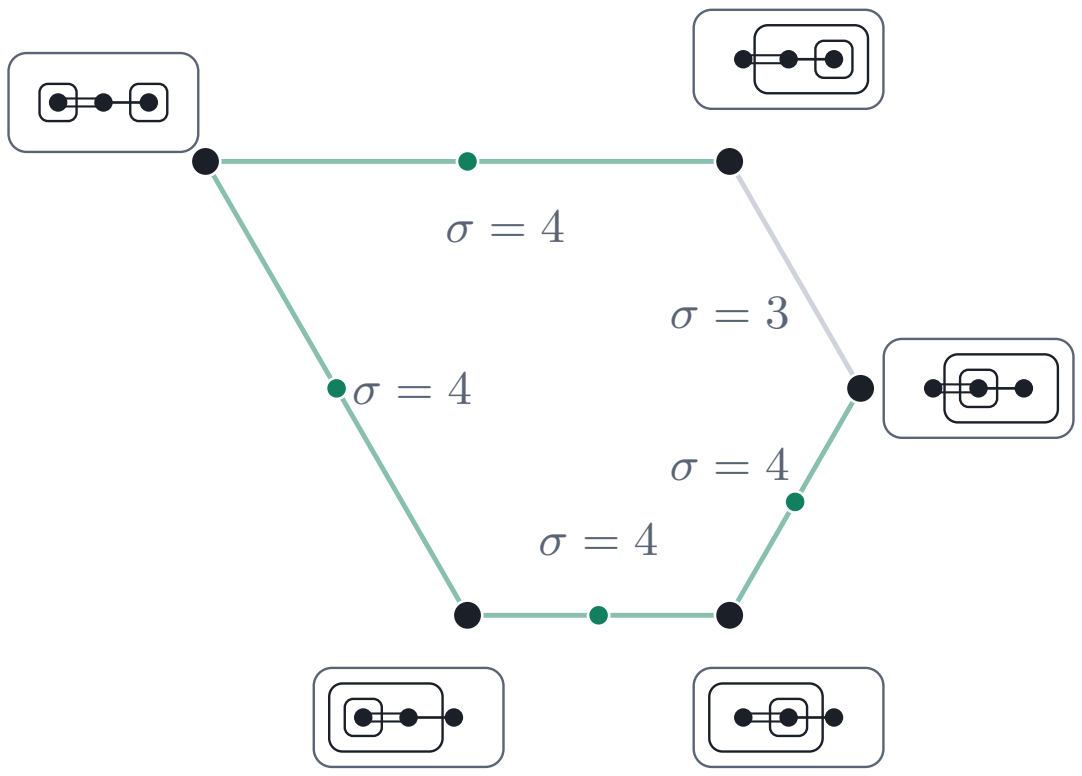

\begin{example}
\label{ex:B3all}
Let $W$ be our $B_3$ running example from \cref{ex:B3cube}. The cubical tubings of $\Gamma$ are the five maximal tubings of \cref{ex:B3nested} and the five edges $\{S, 12\}$, $\{S, 23\}$, $\{S, 1\}$, $\{S, 2\}$ and $\{S, 3\}$ of the pentagon. For the first two edges the two-element shell lies in $T = 12$ ($\sigma(12, s_1, s_2) = 4$) and $T = 23$ ($\sigma(23, s_2, s_3) = 3$) respectively. For the other three edges it lies in $T = S$ and $J$ is a single node. The case $J = \{s_1\}$ was worked out in \cref{ex:B3flat}, and we got $\sigma(S, s_2, s_3) = 4$. For $J = \{s_2\}$ we have (ratios are $0, 1/2, 1, \infty$)
\begin{align*}
L_1 &= \{\alpha_1,\ \alpha_1+\alpha_2,\ 2\alpha_1+\alpha_2\}, \\
L_2 &= \{2\alpha_1+\alpha_2+\alpha_3,\ 2\alpha_1+2\alpha_2+\alpha_3\}, \\
L_3 &= \{\alpha_1+\alpha_2+\alpha_3\}, \\
L_4 &= \{\alpha_3,\ \alpha_2+\alpha_3\},
\end{align*}
so we get $\sigma(S, s_1, s_3) = 4$. For $J = \{s_3\}$ we have (ratios are again $0, 1/2, 1, \infty$)
\begin{align*}
L_1 &= \{\alpha_1\}, \\
L_2 &= \{2\alpha_1+\alpha_2,\ 2\alpha_1+\alpha_2+\alpha_3\}, \\
L_3 &= \{\alpha_1+\alpha_2,\ \alpha_1+\alpha_2+\alpha_3,\ 2\alpha_1+2\alpha_2+\alpha_3\}, \\
L_4 &= \{\alpha_2,\ \alpha_2+\alpha_3\},
\end{align*}
so we get $\sigma(S, s_1, s_2) = 4$. By \cref{thm:cubicaltubing}, the edge $\{S, 23\}$ carries no based cube, and each of the other four edges carries exactly one. In each case the block containing $\alpha_p$ is $L_1 \sqcup L_2$ and the block containing $\alpha_q$ is $L_3 \sqcup L_4$. Writing the block containing $\alpha_p$ first, we get
\begin{align*}
\{S, 1\} &\colon \{\alpha_2,\ \alpha_1+\alpha_2,\ 2\alpha_1+\alpha_2,\ 2\alpha_1+2\alpha_2+\alpha_3\} \\
&\qquad \sqcup \{\alpha_3,\ \alpha_2+\alpha_3,\ \alpha_1+\alpha_2+\alpha_3,\ 2\alpha_1+\alpha_2+\alpha_3\} \sqcup \{\alpha_1\}, \\
\{S, 2\} &\colon \{\alpha_1,\ \alpha_1+\alpha_2,\ 2\alpha_1+\alpha_2,\ 2\alpha_1+\alpha_2+\alpha_3,\ 2\alpha_1+2\alpha_2+\alpha_3\} \\
&\qquad \sqcup \{\alpha_3,\ \alpha_2+\alpha_3,\ \alpha_1+\alpha_2+\alpha_3\} \sqcup \{\alpha_2\}, \\
\{S, 3\} &\colon \{\alpha_1,\ 2\alpha_1+\alpha_2,\ 2\alpha_1+\alpha_2+\alpha_3\} \\
&\qquad \sqcup \{\alpha_2,\ \alpha_1+\alpha_2,\ \alpha_2+\alpha_3,\ \alpha_1+\alpha_2+\alpha_3,\ 2\alpha_1+2\alpha_2+\alpha_3\} \sqcup \{\alpha_3\}, \\
\{S, 12\} &\colon \{\alpha_1,\ 2\alpha_1+\alpha_2\} \sqcup \{\alpha_2,\ \alpha_1+\alpha_2\} \sqcup B(S,3).
\end{align*}
The first one is the cube $\C$ of \cref{ex:B3cube}. In total, $W$ has $5 + 4 = 9$ based cubes, placed on $K(\Gamma)$ as in \cref{fig:B3twisted}.
\end{example}

\begin{example}
\label{ex:H4}
Let $W$ be of type $H_4$, with $S = \{s_1, s_2, s_3, s_4\}$, $m_{12} = 5$ and $m_{23} = m_{34} = 3$. Since $\Gamma$ is a path on four nodes, $K(\Gamma)$ is the three-dimensional associahedron. It has $14$ vertices, $21$ edges and $9$ two-dimensional faces ($6$ pentagons and $3$ squares). The squares are the cubical tubings $\{S, 12\}$, $\{S, 23\}$ and $\{S, 34\}$. The product structure of \cref{thm:cubicaltubing}(3) can be checked on the square $\{S, 12\}$. Its $14$ cubes form a $7 \times 2$ grid, the index $k_S$ is shared with the $7$ cubes on each of the edges $\{S, 12, 1\}$ and $\{S, 12, 2\}$, and the index $k_{12}$ is shared with the $2$ cubes on each of the edges $\{S, 12, 123\}$ and $\{S, 12, 4\}$.
\end{example}

\subsection{Loosely separable elements}

We define loosely separable elements, dropping the parabolic requirement in the description of \cref{thm:nestedvertices}(2).

\begin{definition}
A \emph{cube decomposition} of an inversion set $\Phi_w$ is a partition of $\Phi_w$ into inversion sets each containing exactly one simple root. We call $w \in W$ \emph{loosely separable} if both $\Phi_w$ and $\Phi^+ \setminus \Phi_w = \Phi_{ww_0}$ (from \eqref{eqn:complement}) admit cube decompositions.
\end{definition}

In the terminology of \cite{Hohlweg-Pons26}, a cube decomposition of $\Phi_w$ is a proper partition of $w$ with the maximal number of blocks. From \cref{thm:nestedvertices}(2), an element $w$ is separable if and only if both $\Phi_w$ and $\Phi^+ \setminus \Phi_w$ have cube decompositions with every block parabolic. It turns out that this loosely separable condition describes exactly the set of all vertices that can appear on based cubes.

\begin{theorem}\label{thm:vertices}
An element $w \in W$ is a vertex of some based cube if and only if $w$ is loosely separable.
\end{theorem}
\begin{proof}
If $w$ is a vertex of a based cube, then $\Phi_w$ is a union of blocks of the cube by \cref{lem:basedcube}. These blocks give a cube decomposition of $\Phi_w$ and the remaining blocks give a cube decomposition of $\Phi^+ \setminus \Phi_w$.

Conversely, let $\Phi_w = B_1 \sqcup \cdots \sqcup B_k$ and $\Phi^+ \setminus \Phi_w = B_{k+1} \sqcup \cdots \sqcup B_r$ be cube decompositions. Then $\Phi^+ = B_1 \sqcup \cdots \sqcup B_r$ is a partition into inversion sets each containing exactly one simple root, and every union of blocks is an inversion set by \cref{lem:subunion}. Hence the blocks form a based cube by \cref{lem:basedcube}, and $w$ is its vertex corresponding to $B_1 \sqcup \cdots \sqcup B_k$.
\end{proof}

\section{Further questions}
\label{sec:questions}

\cref{thm:structure} is proved for Weyl groups by the quotient root system machinery of \cite{Dimitrov-etal25} and for $I_2(m)$ by hand, and is computer checked for $H_3$ and $H_4$ \cite{GO26code}.

\begin{question}
Is there a uniform proof of \cref{thm:structure} for all finite Coxeter groups?
\end{question}

The number of based cubes is known in every type. The generating functions in types $A$, $B$, $C$ and $D$ and the counts for the exceptional Weyl groups are computed in \cite[\S 10.1]{Dimitrov-etal25}. There are $m - 1$ based cubes for $I_2(m)$, and $16$ and $118$ for $H_3$ and $H_4$ \cite{GO26code}. \cref{cor:count} refines these counts by the cubical face on which a based cube lies.

\begin{question}
Can the generating functions of \cite{Dimitrov-etal25} be derived from the sum over cubical tubings in \cref{cor:count}?
\end{question}

For infinite $W$ there are no Coxeter cubes of dimension $r$ \cite{Gimenez25}, but partitions of $\Phi_w$ into inversion sets make sense for every $w$ \cite{Hohlweg-Pons26}. The replacement for inversion sets we can use is the class of biclosed sets, which coincide with inversion sets for finite $W$ \cite[Lemma~4.1]{Dyer19}, \cite[\S 2]{Hohlweg-Labbe16}. Barkley and Speyer gave combinatorial models in affine type \cite{Barkley-Speyer24}.

\begin{question}
Is there an analogue of \cref{thm:cubicaltubing} in affine type using biclosed sets?
\end{question}

By \cite{Ardila-Reiner-Williams06} the graph associahedron $K(\Gamma)$ is dual to the positive Bergman complex of the Coxeter arrangement $\A$. The faces are indexed by the flats coming from the tubes of the Coxeter diagram. It would be interesting to check if the remaining flats could be expressed using our idea of expressing non-nested cubes as points.

\begin{question}
Is there a description of the point configuration of based cubes in terms of the Bergman complex of the Coxeter arrangement?
\end{question}

\section*{Acknowledgements}
We thank Matthew Dyer for pointing us to \cite{Gaetz-Gao20split} and \cite{Dimitrov-etal25}. The latter led us to \cite{Hohlweg-Pons26}, and \cref{lem:subunion} from these papers is used repeatedly throughout the paper. The second author thanks Edward Richmond and Nathan Lesnevich for discussions on Richardson varieties, where he learned of the results and ideas of \cite{Dyer90} used here.

\section*{Use of AI tools}
The identification of based cubes with points of the graph associahedron, the idea and notion of loosely separable elements, and the main proof strategy, in particular the connection with the Coxeter arrangement, were developed by the authors. Claude (Anthropic, mainly the models Claude Fable 5.1 and Claude Opus 5.5) was used during the preparation of the manuscript to assist with some technical details, mainly \cref{subsec:ranktwo} (the arguments were produced in root-theoretic terms and then completely rewritten by the authors in the language of the Coxeter arrangement), the computations behind the examples and the computer verifications in \cite{GO26code}, and with the preparation of the figures. It was also used to help check and polish some arguments. All AI-assisted arguments were independently verified and revised by the authors. The authors take full responsibility for all content in this paper.
\printbibliography
\end{document}

%% file: fig/cube-B3.tex
\definecolor{tokfg}{HTML}{1B1F27}
\definecolor{tokmuted}{HTML}{5B6474}
\definecolor{tokline}{HTML}{CFD4DC}
\definecolor{tokfaint}{HTML}{E3E6EC}
\definecolor{toksoft}{HTML}{DDE3EC}
\definecolor{tokbx}{HTML}{2A6FC9}
\definecolor{tokby}{HTML}{D4582A}
\definecolor{toktw}{HTML}{12805F}
\definecolor{toktwsoft}{HTML}{D7EFE6}

\tikzset{
  kedge/.style={draw=tokline,line width=0.92pt,line cap=round},
  kedgetw/.style={draw=toktw,opacity=.5,line width=0.92pt,line cap=round},
  vdot/.style={circle,fill=tokfg,draw=white,line width=.5pt,inner sep=0pt,minimum size=0.195cm},
  tdot/.style={circle,fill=toktw,draw=white,line width=.5pt,inner sep=0pt,minimum size=0.140cm},
  ttdot/.style={rectangle,rounded corners=1pt,fill=toktw,draw=white,line width=.5pt,inner sep=0pt,minimum size=0.146cm},
  sq/.style={fill=toktwsoft,draw=none},
  lab/.style={font=\small,text=tokfg},
  mlab/.style={font=\footnotesize,text=tokmuted},
}

\begin{tikzpicture}[>=stealth]
\draw[tokfg,line width=.9pt,->,shorten >=11pt,shorten <=11pt] (-2.60,1.60) -- (0.00,0.00);
\draw[tokbx,line width=.9pt,->,shorten >=11pt,shorten <=11pt] (0.00,1.60) -- (0.00,0.00);
\draw[tokby,line width=.9pt,->,shorten >=11pt,shorten <=11pt] (2.60,1.60) -- (0.00,0.00);
\draw[tokbx,line width=.9pt,->,shorten >=11pt,shorten <=11pt] (-2.60,3.20) -- (-2.60,1.60);
\draw[tokby,line width=.9pt,->,shorten >=11pt,shorten <=11pt] (0.00,3.20) -- (-2.60,1.60);
\draw[tokfg,line width=.9pt,->,shorten >=11pt,shorten <=11pt] (-2.60,3.20) -- (0.00,1.60);
\draw[tokby,line width=.9pt,->,shorten >=11pt,shorten <=11pt] (2.60,3.20) -- (0.00,1.60);
\draw[tokfg,line width=.9pt,->,shorten >=11pt,shorten <=11pt] (0.00,3.20) -- (2.60,1.60);
\draw[tokbx,line width=.9pt,->,shorten >=11pt,shorten <=11pt] (2.60,3.20) -- (2.60,1.60);
\draw[tokby,line width=.9pt,->,shorten >=11pt,shorten <=11pt] (0.00,4.80) -- (-2.60,3.20);
\draw[tokbx,line width=.9pt,->,shorten >=11pt,shorten <=11pt] (0.00,4.80) -- (0.00,3.20);
\draw[tokfg,line width=.9pt,->,shorten >=11pt,shorten <=11pt] (0.00,4.80) -- (2.60,3.20);
\node[lab,fill=white,inner sep=1.5pt,rounded corners=1pt] at (0.00,0.00) {$e$};
\node[lab,fill=white,inner sep=1.5pt,rounded corners=1pt] at (-2.60,1.60) {$s_1$};
\node[lab,fill=white,inner sep=1.5pt,rounded corners=1pt] at (0.00,1.60) {$s_2s_1s_2s_3$};
\node[lab,fill=white,inner sep=1.5pt,rounded corners=1pt] at (2.60,1.60) {$s_3s_2s_1s_2$};
\node[lab,fill=white,inner sep=1.5pt,rounded corners=1pt] at (-2.60,3.20) {$s_1s_2s_1s_2s_3$};
\node[lab,fill=white,inner sep=1.5pt,rounded corners=1pt] at (0.00,3.20) {$s_1s_3s_2s_1s_2$};
\node[lab,fill=white,inner sep=1.5pt,rounded corners=1pt] at (2.60,3.20) {$s_2s_1s_2s_3s_2s_1s_2s_3$};
\node[lab,fill=white,inner sep=1.5pt,rounded corners=1pt] at (0.00,4.80) {$w_0 = s_1s_2s_1s_2s_3s_2s_1s_2s_3$};
\node[mlab,text=tokfg] at (-1.62,0.62) {$x_1$};
\node[mlab,text=tokbx] at (0.28,0.80) {$x_2$};
\node[mlab,text=tokby] at (1.62,0.62) {$x_3$};
\node[mlab,anchor=west] at (3.9,4.3) {edges of the parallel class of};
\draw[tokfg,line width=.9pt,->] (3.95,3.90) -- (4.55,3.90);
\node[mlab,anchor=west] at (4.6,3.90) {$B_1$};
\draw[tokbx,line width=.9pt,->] (3.95,3.48) -- (4.55,3.48);
\node[mlab,anchor=west] at (4.6,3.48) {$B_2$};
\draw[tokby,line width=.9pt,->] (3.95,3.06) -- (4.55,3.06);
\node[mlab,anchor=west] at (4.6,3.06) {$B_3$};
\end{tikzpicture}

%% file: fig/pentagon-B3-nested.tex
\definecolor{tokfg}{HTML}{1B1F27}
\definecolor{tokmuted}{HTML}{5B6474}
\definecolor{tokline}{HTML}{CFD4DC}
\definecolor{tokfaint}{HTML}{E3E6EC}
\definecolor{toksoft}{HTML}{DDE3EC}
\definecolor{tokbx}{HTML}{2A6FC9}
\definecolor{tokby}{HTML}{D4582A}
\definecolor{toktw}{HTML}{12805F}
\definecolor{toktwsoft}{HTML}{D7EFE6}

\tikzset{
  kedge/.style={draw=tokline,line width=0.92pt,line cap=round},
  kedgetw/.style={draw=toktw,opacity=.5,line width=0.92pt,line cap=round},
  vdot/.style={circle,fill=tokfg,draw=white,line width=.5pt,inner sep=0pt,minimum size=0.195cm},
  tdot/.style={circle,fill=toktw,draw=white,line width=.5pt,inner sep=0pt,minimum size=0.140cm},
  ttdot/.style={rectangle,rounded corners=1pt,fill=toktw,draw=white,line width=.5pt,inner sep=0pt,minimum size=0.146cm},
  sq/.style={fill=toktwsoft,draw=none},
  lab/.style={font=\small,text=tokfg},
  mlab/.style={font=\footnotesize,text=tokmuted},
}

\begin{tikzpicture}
\draw[kedge] (-2.600,1.501) -- (-0.867,-1.501);
\draw[kedge] (-2.600,1.501) -- (0.867,1.501);
\draw[kedge] (-0.867,-1.501) -- (0.867,-1.501);
\draw[kedge] (0.867,-1.501) -- (1.733,0.000);
\draw[kedge] (1.733,0.000) -- (0.867,1.501);
\node[vdot] at (-2.600,1.501) {};
\draw[tokfg,line width=.4pt,double,double distance=1.1pt] (-3.575,1.891) -- (-3.275,1.891);
\draw[tokfg,line width=.5pt] (-3.275,1.891) -- (-2.975,1.891);
\fill[tokfg] (-3.575,1.891) circle (0.061);
\fill[tokfg] (-3.275,1.891) circle (0.061);
\fill[tokfg] (-2.975,1.891) circle (0.061);
\draw[tokfg,line width=.45pt,rounded corners=1.6pt] (-3.098,1.768) rectangle (-2.853,2.014);
\draw[tokfg,line width=.45pt,rounded corners=1.6pt] (-3.698,1.768) rectangle (-3.453,2.014);
\draw[tokmuted,line width=.45pt,rounded corners=3.4pt] (-3.903,1.564) rectangle (-2.648,2.218);
\node[vdot] at (-0.867,-1.501) {};
\draw[tokfg,line width=.4pt,double,double distance=1.1pt] (-1.557,-2.177) -- (-1.257,-2.177);
\draw[tokfg,line width=.5pt] (-1.257,-2.177) -- (-0.957,-2.177);
\fill[tokfg] (-1.557,-2.177) circle (0.061);
\fill[tokfg] (-1.257,-2.177) circle (0.061);
\fill[tokfg] (-0.957,-2.177) circle (0.061);
\draw[tokfg,line width=.45pt,rounded corners=1.6pt] (-1.679,-2.299) rectangle (-1.434,-2.054);
\draw[tokfg,line width=.45pt,rounded corners=2.5pt] (-1.782,-2.402) rectangle (-1.032,-1.952);
\draw[tokmuted,line width=.45pt,rounded corners=3.4pt] (-1.884,-2.504) rectangle (-0.629,-1.849);
\node[vdot] at (0.867,-1.501) {};
\draw[tokfg,line width=.4pt,double,double distance=1.1pt] (0.957,-2.177) -- (1.257,-2.177);
\draw[tokfg,line width=.5pt] (1.257,-2.177) -- (1.557,-2.177);
\fill[tokfg] (0.957,-2.177) circle (0.061);
\fill[tokfg] (1.257,-2.177) circle (0.061);
\fill[tokfg] (1.557,-2.177) circle (0.061);
\draw[tokfg,line width=.45pt,rounded corners=1.6pt] (1.134,-2.299) rectangle (1.379,-2.054);
\draw[tokfg,line width=.45pt,rounded corners=2.5pt] (0.732,-2.402) rectangle (1.482,-1.952);
\draw[tokmuted,line width=.45pt,rounded corners=3.4pt] (0.629,-2.504) rectangle (1.884,-1.849);
\node[vdot] at (1.733,0.000) {};
\draw[tokfg,line width=.4pt,double,double distance=1.1pt] (2.213,0.000) -- (2.513,0.000);
\draw[tokfg,line width=.5pt] (2.513,0.000) -- (2.813,0.000);
\fill[tokfg] (2.213,0.000) circle (0.061);
\fill[tokfg] (2.513,0.000) circle (0.061);
\fill[tokfg] (2.813,0.000) circle (0.061);
\draw[tokfg,line width=.45pt,rounded corners=1.6pt] (2.391,-0.123) rectangle (2.636,0.123);
\draw[tokfg,line width=.45pt,rounded corners=2.5pt] (2.288,-0.225) rectangle (3.038,0.225);
\draw[tokmuted,line width=.45pt,rounded corners=3.4pt] (1.886,-0.327) rectangle (3.141,0.327);
\node[vdot] at (0.867,1.501) {};
\draw[tokfg,line width=.4pt,double,double distance=1.1pt] (0.957,2.177) -- (1.257,2.177);
\draw[tokfg,line width=.5pt] (1.257,2.177) -- (1.557,2.177);
\fill[tokfg] (0.957,2.177) circle (0.061);
\fill[tokfg] (1.257,2.177) circle (0.061);
\fill[tokfg] (1.557,2.177) circle (0.061);
\draw[tokfg,line width=.45pt,rounded corners=1.6pt] (1.434,2.054) rectangle (1.679,2.299);
\draw[tokfg,line width=.45pt,rounded corners=2.5pt] (1.032,1.952) rectangle (1.782,2.402);
\draw[tokmuted,line width=.45pt,rounded corners=3.4pt] (0.629,1.849) rectangle (1.884,2.504);
\end{tikzpicture}

%% file: fig/flat-B3.tex
\definecolor{tokfg}{HTML}{1B1F27}
\definecolor{tokmuted}{HTML}{5B6474}
\definecolor{tokline}{HTML}{CFD4DC}
\definecolor{tokfaint}{HTML}{E3E6EC}
\definecolor{toksoft}{HTML}{DDE3EC}
\definecolor{tokbx}{HTML}{2A6FC9}
\definecolor{tokby}{HTML}{D4582A}
\definecolor{toktw}{HTML}{12805F}
\definecolor{toktwsoft}{HTML}{D7EFE6}

\tikzset{
  kedge/.style={draw=tokline,line width=0.92pt,line cap=round},
  kedgetw/.style={draw=toktw,opacity=.5,line width=0.92pt,line cap=round},
  vdot/.style={circle,fill=tokfg,draw=white,line width=.5pt,inner sep=0pt,minimum size=0.195cm},
  tdot/.style={circle,fill=toktw,draw=white,line width=.5pt,inner sep=0pt,minimum size=0.140cm},
  ttdot/.style={rectangle,rounded corners=1pt,fill=toktw,draw=white,line width=.5pt,inner sep=0pt,minimum size=0.146cm},
  sq/.style={fill=toktwsoft,draw=none},
  lab/.style={font=\small,text=tokfg},
  mlab/.style={font=\footnotesize,text=tokmuted},
}

\begin{tikzpicture}
\fill[toktwsoft] (0,0) -- (135:2.55) arc (135:180:2.55) -- cycle;
\fill[toktwsoft] (0,0) -- (315:2.55) arc (315:360:2.55) -- cycle;
\draw[tokbx,line width=.9pt] (90:2.60) -- (270:2.60);
\node[lab,text=tokbx] at (90:3.10) {$\ell_1$};
\draw[tokfg,line width=.9pt] (135:2.60) -- (315:2.60);
\node[lab,text=tokfg] at (135:3.10) {$\ell_2$};
\draw[tokfg,line width=.9pt] (180:2.60) -- (360:2.60);
\node[lab,text=tokfg] at (180:3.10) {$\ell_3$};
\draw[tokby,line width=.9pt] (225:2.60) -- (405:2.60);
\node[lab,text=tokby] at (225:3.10) {$\ell_4$};
\node[mlab,text=toktw] at (157.5:1.87) {$B_2$};
\node[mlab,text=toktw] at (337.5:1.87) {$B_3$};
\end{tikzpicture}

%% file: fig/hemisphere-B3.tex
\definecolor{tokfg}{HTML}{1B1F27}
\definecolor{tokmuted}{HTML}{5B6474}
\definecolor{tokline}{HTML}{CFD4DC}
\definecolor{tokfaint}{HTML}{E3E6EC}
\definecolor{toksoft}{HTML}{DDE3EC}
\definecolor{tokbx}{HTML}{2A6FC9}
\definecolor{tokby}{HTML}{D4582A}
\definecolor{toktw}{HTML}{12805F}
\definecolor{toktwsoft}{HTML}{D7EFE6}

\tikzset{
  kedge/.style={draw=tokline,line width=0.92pt,line cap=round},
  kedgetw/.style={draw=toktw,opacity=.5,line width=0.92pt,line cap=round},
  vdot/.style={circle,fill=tokfg,draw=white,line width=.5pt,inner sep=0pt,minimum size=0.195cm},
  tdot/.style={circle,fill=toktw,draw=white,line width=.5pt,inner sep=0pt,minimum size=0.140cm},
  ttdot/.style={rectangle,rounded corners=1pt,fill=toktw,draw=white,line width=.5pt,inner sep=0pt,minimum size=0.146cm},
  sq/.style={fill=toktwsoft,draw=none},
  lab/.style={font=\small,text=tokfg},
  mlab/.style={font=\footnotesize,text=tokmuted},
}

\begin{tikzpicture}
\fill[fill=toksoft] (1.501,1.501) -- (1.518,1.518) -- (1.534,1.534) -- (1.550,1.550) -- (1.565,1.565) -- (1.580,1.580) -- (1.595,1.595) -- (1.609,1.609) -- (1.623,1.623) -- (1.636,1.636) -- (1.649,1.649) -- (1.662,1.662) -- (1.674,1.674) -- (1.686,1.686) -- (1.697,1.697) -- (1.708,1.708) -- (1.719,1.719) -- (1.729,1.729) -- (1.738,1.738) -- (1.748,1.748) -- (1.756,1.756) -- (1.765,1.765) -- (1.773,1.773) -- (1.780,1.780) -- (1.787,1.787) -- (1.794,1.794) -- (1.800,1.800) -- (1.806,1.806) -- (1.811,1.811) -- (1.816,1.816) -- (1.820,1.820) -- (1.824,1.824) -- (1.827,1.827) -- (1.830,1.830) -- (1.833,1.833) -- (1.835,1.835) -- (1.836,1.836) -- (1.838,1.838) -- (1.838,1.838) -- (1.838,1.838) -- (1.838,1.838) -- (1.801,1.875) -- (1.763,1.911) -- (1.724,1.946) -- (1.685,1.980) -- (1.644,2.014) -- (1.603,2.047) -- (1.562,2.079) -- (1.520,2.110) -- (1.477,2.140) -- (1.434,2.169) -- (1.390,2.197) -- (1.345,2.225) -- (1.300,2.252) -- (1.254,2.277) -- (1.208,2.302) -- (1.162,2.326) -- (1.115,2.349) -- (1.067,2.371) -- (1.019,2.392) -- (0.971,2.412) -- (0.922,2.431) -- (0.873,2.449) -- (0.823,2.466) -- (0.774,2.482) -- (0.723,2.497) -- (0.673,2.511) -- (0.622,2.524) -- (0.571,2.536) -- (0.520,2.547) -- (0.469,2.557) -- (0.417,2.566) -- (0.365,2.574) -- (0.313,2.581) -- (0.261,2.587) -- (0.209,2.592) -- (0.157,2.595) -- (0.105,2.598) -- (0.052,2.599) -- (0.000,2.600) -- (0.000,2.600) -- (0.045,2.599) -- (0.090,2.597) -- (0.135,2.593) -- (0.180,2.588) -- (0.225,2.581) -- (0.269,2.572) -- (0.314,2.562) -- (0.358,2.550) -- (0.402,2.537) -- (0.446,2.522) -- (0.489,2.506) -- (0.533,2.488) -- (0.576,2.469) -- (0.618,2.449) -- (0.660,2.426) -- (0.702,2.403) -- (0.744,2.378) -- (0.785,2.351) -- (0.825,2.323) -- (0.865,2.294) -- (0.905,2.264) -- (0.943,2.232) -- (0.982,2.198) -- (1.020,2.163) -- (1.057,2.128) -- (1.093,2.090) -- (1.129,2.052) -- (1.164,2.012) -- (1.199,1.971) -- (1.233,1.929) -- (1.266,1.886) -- (1.298,1.841) -- (1.329,1.796) -- (1.360,1.749) -- (1.390,1.702) -- (1.419,1.653) -- (1.447,1.603) -- (1.475,1.553) -- (1.501,1.501) -- cycle;
\fill[fill=toktwsoft] (-1.501,1.501) -- (-1.518,1.518) -- (-1.534,1.534) -- (-1.550,1.550) -- (-1.565,1.565) -- (-1.580,1.580) -- (-1.595,1.595) -- (-1.609,1.609) -- (-1.623,1.623) -- (-1.636,1.636) -- (-1.649,1.649) -- (-1.662,1.662) -- (-1.674,1.674) -- (-1.686,1.686) -- (-1.697,1.697) -- (-1.708,1.708) -- (-1.719,1.719) -- (-1.729,1.729) -- (-1.738,1.738) -- (-1.748,1.748) -- (-1.756,1.756) -- (-1.765,1.765) -- (-1.773,1.773) -- (-1.780,1.780) -- (-1.787,1.787) -- (-1.794,1.794) -- (-1.800,1.800) -- (-1.806,1.806) -- (-1.811,1.811) -- (-1.816,1.816) -- (-1.820,1.820) -- (-1.824,1.824) -- (-1.827,1.827) -- (-1.830,1.830) -- (-1.833,1.833) -- (-1.835,1.835) -- (-1.836,1.836) -- (-1.838,1.838) -- (-1.838,1.838) -- (-1.838,1.838) -- (-1.838,1.838) -- (-1.875,1.801) -- (-1.911,1.763) -- (-1.946,1.724) -- (-1.980,1.685) -- (-2.014,1.644) -- (-2.047,1.603) -- (-2.079,1.562) -- (-2.110,1.520) -- (-2.140,1.477) -- (-2.169,1.434) -- (-2.197,1.390) -- (-2.225,1.345) -- (-2.252,1.300) -- (-2.277,1.254) -- (-2.302,1.208) -- (-2.326,1.162) -- (-2.349,1.115) -- (-2.371,1.067) -- (-2.392,1.019) -- (-2.412,0.971) -- (-2.431,0.922) -- (-2.449,0.873) -- (-2.466,0.823) -- (-2.482,0.774) -- (-2.497,0.723) -- (-2.511,0.673) -- (-2.524,0.622) -- (-2.536,0.571) -- (-2.547,0.520) -- (-2.557,0.469) -- (-2.566,0.417) -- (-2.574,0.365) -- (-2.581,0.313) -- (-2.587,0.261) -- (-2.592,0.209) -- (-2.595,0.157) -- (-2.598,0.105) -- (-2.599,0.052) -- (-2.600,0.000) -- (-2.600,0.000) -- (-2.599,0.045) -- (-2.597,0.090) -- (-2.593,0.135) -- (-2.588,0.180) -- (-2.581,0.225) -- (-2.572,0.269) -- (-2.562,0.314) -- (-2.550,0.358) -- (-2.537,0.402) -- (-2.522,0.446) -- (-2.506,0.489) -- (-2.488,0.533) -- (-2.469,0.576) -- (-2.449,0.618) -- (-2.426,0.660) -- (-2.403,0.702) -- (-2.378,0.744) -- (-2.351,0.785) -- (-2.323,0.825) -- (-2.294,0.865) -- (-2.264,0.905) -- (-2.232,0.943) -- (-2.198,0.982) -- (-2.163,1.020) -- (-2.128,1.057) -- (-2.090,1.093) -- (-2.052,1.129) -- (-2.012,1.164) -- (-1.971,1.199) -- (-1.929,1.233) -- (-1.886,1.266) -- (-1.841,1.298) -- (-1.796,1.329) -- (-1.749,1.360) -- (-1.702,1.390) -- (-1.653,1.419) -- (-1.603,1.447) -- (-1.553,1.475) -- (-1.501,1.501) -- cycle;
\draw[draw=tokby,line width=.9pt] (0.000,2.600) -- (-0.000,2.596) -- (-0.000,2.586) -- (-0.000,2.568) -- (-0.000,2.543) -- (-0.000,2.511) -- (-0.000,2.473) -- (-0.000,2.427) -- (-0.000,2.375) -- (-0.000,2.317) -- (-0.000,2.252) -- (-0.000,2.181) -- (-0.000,2.103) -- (-0.000,2.021) -- (-0.000,1.932) -- (-0.000,1.838) -- (-0.000,1.740) -- (-0.000,1.636) -- (-0.000,1.528) -- (-0.000,1.416) -- (-0.000,1.300) -- (-0.000,1.180) -- (-0.000,1.058) -- (-0.000,0.932) -- (-0.000,0.803) -- (-0.000,0.673) -- (-0.000,0.541) -- (-0.000,0.407) -- (-0.000,0.272) -- (-0.000,0.136) -- (-0.000,0.000) -- (0.000,-0.136) -- (0.000,-0.272) -- (0.000,-0.407) -- (0.000,-0.541) -- (0.000,-0.673) -- (0.000,-0.803) -- (0.000,-0.932) -- (0.000,-1.058) -- (0.000,-1.180) -- (0.000,-1.300) -- (0.000,-1.416) -- (0.000,-1.528) -- (0.000,-1.636) -- (0.000,-1.740) -- (0.000,-1.838) -- (0.000,-1.932) -- (0.000,-2.021) -- (0.000,-2.103) -- (0.000,-2.181) -- (0.000,-2.252) -- (0.000,-2.317) -- (0.000,-2.375) -- (0.000,-2.427) -- (0.000,-2.473) -- (0.000,-2.511) -- (0.000,-2.543) -- (0.000,-2.568) -- (0.000,-2.586) -- (0.000,-2.596) -- (0.000,-2.600) -- (0.000,-2.600);
\draw[draw=tokline,line width=.6pt] (-2.600,0.000) -- (-2.596,-0.000) -- (-2.586,-0.000) -- (-2.568,-0.000) -- (-2.543,-0.000) -- (-2.511,-0.000) -- (-2.473,-0.000) -- (-2.427,-0.000) -- (-2.375,-0.000) -- (-2.317,-0.000) -- (-2.252,-0.000) -- (-2.181,-0.000) -- (-2.103,-0.000) -- (-2.021,-0.000) -- (-1.932,-0.000) -- (-1.838,-0.000) -- (-1.740,-0.000) -- (-1.636,-0.000) -- (-1.528,-0.000) -- (-1.416,-0.000) -- (-1.300,-0.000) -- (-1.180,-0.000) -- (-1.058,-0.000) -- (-0.932,-0.000) -- (-0.803,-0.000) -- (-0.673,-0.000) -- (-0.541,-0.000) -- (-0.407,-0.000) -- (-0.272,-0.000) -- (-0.136,-0.000) -- (-0.000,-0.000) -- (0.136,0.000) -- (0.272,0.000) -- (0.407,0.000) -- (0.541,0.000) -- (0.673,0.000) -- (0.803,0.000) -- (0.932,0.000) -- (1.058,0.000) -- (1.180,0.000) -- (1.300,0.000) -- (1.416,0.000) -- (1.528,0.000) -- (1.636,0.000) -- (1.740,0.000) -- (1.838,0.000) -- (1.932,0.000) -- (2.021,0.000) -- (2.103,0.000) -- (2.181,0.000) -- (2.252,0.000) -- (2.317,0.000) -- (2.375,0.000) -- (2.427,0.000) -- (2.473,0.000) -- (2.511,0.000) -- (2.543,0.000) -- (2.568,0.000) -- (2.586,0.000) -- (2.596,0.000) -- (2.600,0.000) -- (2.600,0.000);
\draw[draw=tokby,line width=.9pt] (-0.000,2.600) -- (0.096,2.596) -- (0.192,2.586) -- (0.288,2.568) -- (0.382,2.543) -- (0.476,2.511) -- (0.568,2.473) -- (0.659,2.427) -- (0.748,2.375) -- (0.835,2.317) -- (0.919,2.252) -- (1.001,2.181) -- (1.081,2.103) -- (1.157,2.021) -- (1.230,1.932) -- (1.300,1.838) -- (1.366,1.740) -- (1.429,1.636) -- (1.487,1.528) -- (1.542,1.416) -- (1.592,1.300) -- (1.638,1.180) -- (1.680,1.058) -- (1.716,0.932) -- (1.748,0.803) -- (1.776,0.673) -- (1.798,0.541) -- (1.816,0.407) -- (1.828,0.272) -- (1.836,0.136) -- (1.838,0.000) -- (1.836,-0.136) -- (1.828,-0.272) -- (1.816,-0.407) -- (1.798,-0.541) -- (1.776,-0.673) -- (1.748,-0.803) -- (1.716,-0.932) -- (1.680,-1.058) -- (1.638,-1.180) -- (1.592,-1.300) -- (1.542,-1.416) -- (1.487,-1.528) -- (1.429,-1.636) -- (1.366,-1.740) -- (1.300,-1.838) -- (1.230,-1.932) -- (1.157,-2.021) -- (1.081,-2.103) -- (1.001,-2.181) -- (0.919,-2.252) -- (0.835,-2.317) -- (0.748,-2.375) -- (0.659,-2.427) -- (0.568,-2.473) -- (0.476,-2.511) -- (0.382,-2.543) -- (0.288,-2.568) -- (0.192,-2.586) -- (0.096,-2.596) -- (0.000,-2.600) -- (0.000,-2.600);
\draw[draw=tokby,line width=.9pt] (0.000,2.600) -- (-0.096,2.596) -- (-0.192,2.586) -- (-0.288,2.568) -- (-0.382,2.543) -- (-0.476,2.511) -- (-0.568,2.473) -- (-0.659,2.427) -- (-0.748,2.375) -- (-0.835,2.317) -- (-0.919,2.252) -- (-1.001,2.181) -- (-1.081,2.103) -- (-1.157,2.021) -- (-1.230,1.932) -- (-1.300,1.838) -- (-1.366,1.740) -- (-1.429,1.636) -- (-1.487,1.528) -- (-1.542,1.416) -- (-1.592,1.300) -- (-1.638,1.180) -- (-1.680,1.058) -- (-1.716,0.932) -- (-1.748,0.803) -- (-1.776,0.673) -- (-1.798,0.541) -- (-1.816,0.407) -- (-1.828,0.272) -- (-1.836,0.136) -- (-1.838,0.000) -- (-1.836,-0.136) -- (-1.828,-0.272) -- (-1.816,-0.407) -- (-1.798,-0.541) -- (-1.776,-0.673) -- (-1.748,-0.803) -- (-1.716,-0.932) -- (-1.680,-1.058) -- (-1.638,-1.180) -- (-1.592,-1.300) -- (-1.542,-1.416) -- (-1.487,-1.528) -- (-1.429,-1.636) -- (-1.366,-1.740) -- (-1.300,-1.838) -- (-1.230,-1.932) -- (-1.157,-2.021) -- (-1.081,-2.103) -- (-1.001,-2.181) -- (-0.919,-2.252) -- (-0.835,-2.317) -- (-0.748,-2.375) -- (-0.659,-2.427) -- (-0.568,-2.473) -- (-0.476,-2.511) -- (-0.382,-2.543) -- (-0.288,-2.568) -- (-0.192,-2.586) -- (-0.096,-2.596) -- (-0.000,-2.600) -- (0.000,-2.600);
\draw[draw=tokline,line width=.6pt] (-2.600,-0.000) -- (-2.596,0.096) -- (-2.586,0.192) -- (-2.568,0.288) -- (-2.543,0.382) -- (-2.511,0.476) -- (-2.473,0.568) -- (-2.427,0.659) -- (-2.375,0.748) -- (-2.317,0.835) -- (-2.252,0.919) -- (-2.181,1.001) -- (-2.103,1.081) -- (-2.021,1.157) -- (-1.932,1.230) -- (-1.838,1.300) -- (-1.740,1.366) -- (-1.636,1.429) -- (-1.528,1.487) -- (-1.416,1.542) -- (-1.300,1.592) -- (-1.180,1.638) -- (-1.058,1.680) -- (-0.932,1.716) -- (-0.803,1.748) -- (-0.673,1.776) -- (-0.541,1.798) -- (-0.407,1.816) -- (-0.272,1.828) -- (-0.136,1.836) -- (-0.000,1.838) -- (0.136,1.836) -- (0.272,1.828) -- (0.407,1.816) -- (0.541,1.798) -- (0.673,1.776) -- (0.803,1.748) -- (0.932,1.716) -- (1.058,1.680) -- (1.180,1.638) -- (1.300,1.592) -- (1.416,1.542) -- (1.528,1.487) -- (1.636,1.429) -- (1.740,1.366) -- (1.838,1.300) -- (1.932,1.230) -- (2.021,1.157) -- (2.103,1.081) -- (2.181,1.001) -- (2.252,0.919) -- (2.317,0.835) -- (2.375,0.748) -- (2.427,0.659) -- (2.473,0.568) -- (2.511,0.476) -- (2.543,0.382) -- (2.568,0.288) -- (2.586,0.192) -- (2.596,0.096) -- (2.600,0.000) -- (2.600,-0.000);
\draw[draw=tokline,line width=.6pt] (-2.600,0.000) -- (-2.596,-0.096) -- (-2.586,-0.192) -- (-2.568,-0.288) -- (-2.543,-0.382) -- (-2.511,-0.476) -- (-2.473,-0.568) -- (-2.427,-0.659) -- (-2.375,-0.748) -- (-2.317,-0.835) -- (-2.252,-0.919) -- (-2.181,-1.001) -- (-2.103,-1.081) -- (-2.021,-1.157) -- (-1.932,-1.230) -- (-1.838,-1.300) -- (-1.740,-1.366) -- (-1.636,-1.429) -- (-1.528,-1.487) -- (-1.416,-1.542) -- (-1.300,-1.592) -- (-1.180,-1.638) -- (-1.058,-1.680) -- (-0.932,-1.716) -- (-0.803,-1.748) -- (-0.673,-1.776) -- (-0.541,-1.798) -- (-0.407,-1.816) -- (-0.272,-1.828) -- (-0.136,-1.836) -- (-0.000,-1.838) -- (0.136,-1.836) -- (0.272,-1.828) -- (0.407,-1.816) -- (0.541,-1.798) -- (0.673,-1.776) -- (0.803,-1.748) -- (0.932,-1.716) -- (1.058,-1.680) -- (1.180,-1.638) -- (1.300,-1.592) -- (1.416,-1.542) -- (1.528,-1.487) -- (1.636,-1.429) -- (1.740,-1.366) -- (1.838,-1.300) -- (1.932,-1.230) -- (2.021,-1.157) -- (2.103,-1.081) -- (2.181,-1.001) -- (2.252,-0.919) -- (2.317,-0.835) -- (2.375,-0.748) -- (2.427,-0.659) -- (2.473,-0.568) -- (2.511,-0.476) -- (2.543,-0.382) -- (2.568,-0.288) -- (2.586,-0.192) -- (2.596,-0.096) -- (2.600,-0.000) -- (2.600,0.000);
\draw[draw=tokline,line width=.6pt] (-1.838,-1.838) -- (-1.836,-1.836) -- (-1.828,-1.828) -- (-1.816,-1.816) -- (-1.798,-1.798) -- (-1.776,-1.776) -- (-1.748,-1.748) -- (-1.716,-1.716) -- (-1.680,-1.680) -- (-1.638,-1.638) -- (-1.592,-1.592) -- (-1.542,-1.542) -- (-1.487,-1.487) -- (-1.429,-1.429) -- (-1.366,-1.366) -- (-1.300,-1.300) -- (-1.230,-1.230) -- (-1.157,-1.157) -- (-1.081,-1.081) -- (-1.001,-1.001) -- (-0.919,-0.919) -- (-0.835,-0.835) -- (-0.748,-0.748) -- (-0.659,-0.659) -- (-0.568,-0.568) -- (-0.476,-0.476) -- (-0.382,-0.382) -- (-0.288,-0.288) -- (-0.192,-0.192) -- (-0.096,-0.096) -- (-0.000,-0.000) -- (0.096,0.096) -- (0.192,0.192) -- (0.288,0.288) -- (0.382,0.382) -- (0.476,0.476) -- (0.568,0.568) -- (0.659,0.659) -- (0.748,0.748) -- (0.835,0.835) -- (0.919,0.919) -- (1.001,1.001) -- (1.081,1.081) -- (1.157,1.157) -- (1.230,1.230) -- (1.300,1.300) -- (1.366,1.366) -- (1.429,1.429) -- (1.487,1.487) -- (1.542,1.542) -- (1.592,1.592) -- (1.638,1.638) -- (1.680,1.680) -- (1.716,1.716) -- (1.748,1.748) -- (1.776,1.776) -- (1.798,1.798) -- (1.816,1.816) -- (1.828,1.828) -- (1.836,1.836) -- (1.838,1.838) -- (1.838,1.838);
\draw[draw=tokby,line width=.9pt] (-1.838,1.838) -- (-1.836,1.836) -- (-1.828,1.828) -- (-1.816,1.816) -- (-1.798,1.798) -- (-1.776,1.776) -- (-1.748,1.748) -- (-1.716,1.716) -- (-1.680,1.680) -- (-1.638,1.638) -- (-1.592,1.592) -- (-1.542,1.542) -- (-1.487,1.487) -- (-1.429,1.429) -- (-1.366,1.366) -- (-1.300,1.300) -- (-1.230,1.230) -- (-1.157,1.157) -- (-1.081,1.081) -- (-1.001,1.001) -- (-0.919,0.919) -- (-0.835,0.835) -- (-0.748,0.748) -- (-0.659,0.659) -- (-0.568,0.568) -- (-0.476,0.476) -- (-0.382,0.382) -- (-0.288,0.288) -- (-0.192,0.192) -- (-0.096,0.096) -- (-0.000,0.000) -- (0.096,-0.096) -- (0.192,-0.192) -- (0.288,-0.288) -- (0.382,-0.382) -- (0.476,-0.476) -- (0.568,-0.568) -- (0.659,-0.659) -- (0.748,-0.748) -- (0.835,-0.835) -- (0.919,-0.919) -- (1.001,-1.001) -- (1.081,-1.081) -- (1.157,-1.157) -- (1.230,-1.230) -- (1.300,-1.300) -- (1.366,-1.366) -- (1.429,-1.429) -- (1.487,-1.487) -- (1.542,-1.542) -- (1.592,-1.592) -- (1.638,-1.638) -- (1.680,-1.680) -- (1.716,-1.716) -- (1.748,-1.748) -- (1.776,-1.776) -- (1.798,-1.798) -- (1.816,-1.816) -- (1.828,-1.828) -- (1.836,-1.836) -- (1.838,-1.838) -- (1.838,-1.838);
\draw[tokfg,line width=.9pt] (0,0) circle (2.60);
\draw[toktw,line width=2.2pt] (135:2.60) arc (135:180:2.60);
\fill[tokfg] (90:2.60) circle (0.05);
\node[lab,text=tokbx] at (90:3.00) {$\ell_1$};
\fill[tokfg] (135:2.60) circle (0.05);
\node[lab,text=tokfg] at (135:3.00) {$\ell_2$};
\fill[tokfg] (180:2.60) circle (0.05);
\node[lab,text=tokfg] at (180:3.00) {$\ell_3$};
\fill[tokfg] (225:2.60) circle (0.05);
\node[lab,text=tokby] at (225:3.00) {$\ell_4$};
\fill[tokfg] (270:2.60) circle (0.05);
\fill[tokfg] (315:2.60) circle (0.05);
\fill[tokfg] (0:2.60) circle (0.05);
\fill[tokfg] (45:2.60) circle (0.05);
\node[lab] at (52:2.44) {$C$};
\node[lab,text=toktw,inner sep=0pt] at (140.5:2.41) {$vC$};
\node[mlab,text=toktw] at (157.5:2.95) {$B_2$};
\node[mlab,text=tokby,fill=white,inner sep=1pt] at (2.05,0.95) {$H_{\alpha_2}$};
\node[mlab,text=tokby,fill=white,inner sep=1pt] at (-1.45,-1.8) {$H_{2\alpha_1+\alpha_2}$};
\node[mlab,text=tokby] at (270:2.95) {$H_{\alpha_1+\alpha_2}$};
\node[mlab,text=tokby] at (313:3.35) {$H_{2\alpha_1+2\alpha_2+\alpha_3}$};
\node[mlab] at (352:3.1) {$X_J$};
\end{tikzpicture}

%% file: fig/pentagon-B3-twisted.tex
\definecolor{tokfg}{HTML}{1B1F27}
\definecolor{tokmuted}{HTML}{5B6474}
\definecolor{tokline}{HTML}{CFD4DC}
\definecolor{tokfaint}{HTML}{E3E6EC}
\definecolor{toksoft}{HTML}{DDE3EC}
\definecolor{tokbx}{HTML}{2A6FC9}
\definecolor{tokby}{HTML}{D4582A}
\definecolor{toktw}{HTML}{12805F}
\definecolor{toktwsoft}{HTML}{D7EFE6}

\tikzset{
  kedge/.style={draw=tokline,line width=0.92pt,line cap=round},
  kedgetw/.style={draw=toktw,opacity=.5,line width=0.92pt,line cap=round},
  vdot/.style={circle,fill=tokfg,draw=white,line width=.5pt,inner sep=0pt,minimum size=0.195cm},
  tdot/.style={circle,fill=toktw,draw=white,line width=.5pt,inner sep=0pt,minimum size=0.140cm},
  ttdot/.style={rectangle,rounded corners=1pt,fill=toktw,draw=white,line width=.5pt,inner sep=0pt,minimum size=0.146cm},
  sq/.style={fill=toktwsoft,draw=none},
  lab/.style={font=\small,text=tokfg},
  mlab/.style={font=\footnotesize,text=tokmuted},
}

\begin{tikzpicture}
\draw[kedgetw] (-2.600,1.501) -- (-0.867,-1.501);
\node[mlab] at (-1.233,0.000) {$\sigma=4$};
\draw[kedgetw] (-2.600,1.501) -- (0.867,1.501);
\node[mlab] at (-0.617,1.068) {$\sigma=4$};
\draw[kedgetw] (-0.867,-1.501) -- (0.867,-1.501);
\node[mlab] at (0.000,-1.001) {$\sigma=4$};
\draw[kedgetw] (0.867,-1.501) -- (1.733,0.000);
\node[mlab] at (0.867,-0.501) {$\sigma=4$};
\draw[kedge] (1.733,0.000) -- (0.867,1.501);
\node[mlab] at (0.867,0.501) {$\sigma=3$};
\node[tdot] at (-1.733,0.000) {};
\node[tdot] at (-0.867,1.501) {};
\node[tdot] at (0.000,-1.501) {};
\node[tdot] at (1.300,-0.751) {};
\node[vdot] at (-2.600,1.501) {};
\draw[tokfg,line width=.4pt,double,double distance=1.1pt] (-3.575,1.891) -- (-3.275,1.891);
\draw[tokfg,line width=.5pt] (-3.275,1.891) -- (-2.975,1.891);
\fill[tokfg] (-3.575,1.891) circle (0.061);
\fill[tokfg] (-3.275,1.891) circle (0.061);
\fill[tokfg] (-2.975,1.891) circle (0.061);
\draw[tokfg,line width=.45pt,rounded corners=1.6pt] (-3.098,1.768) rectangle (-2.853,2.014);
\draw[tokfg,line width=.45pt,rounded corners=1.6pt] (-3.698,1.768) rectangle (-3.453,2.014);
\draw[tokmuted,line width=.45pt,rounded corners=3.4pt] (-3.903,1.564) rectangle (-2.648,2.218);
\node[vdot] at (-0.867,-1.501) {};
\draw[tokfg,line width=.4pt,double,double distance=1.1pt] (-1.557,-2.177) -- (-1.257,-2.177);
\draw[tokfg,line width=.5pt] (-1.257,-2.177) -- (-0.957,-2.177);
\fill[tokfg] (-1.557,-2.177) circle (0.061);
\fill[tokfg] (-1.257,-2.177) circle (0.061);
\fill[tokfg] (-0.957,-2.177) circle (0.061);
\draw[tokfg,line width=.45pt,rounded corners=1.6pt] (-1.679,-2.299) rectangle (-1.434,-2.054);
\draw[tokfg,line width=.45pt,rounded corners=2.5pt] (-1.782,-2.402) rectangle (-1.032,-1.952);
\draw[tokmuted,line width=.45pt,rounded corners=3.4pt] (-1.884,-2.504) rectangle (-0.629,-1.849);
\node[vdot] at (0.867,-1.501) {};
\draw[tokfg,line width=.4pt,double,double distance=1.1pt] (0.957,-2.177) -- (1.257,-2.177);
\draw[tokfg,line width=.5pt] (1.257,-2.177) -- (1.557,-2.177);
\fill[tokfg] (0.957,-2.177) circle (0.061);
\fill[tokfg] (1.257,-2.177) circle (0.061);
\fill[tokfg] (1.557,-2.177) circle (0.061);
\draw[tokfg,line width=.45pt,rounded corners=1.6pt] (1.134,-2.299) rectangle (1.379,-2.054);
\draw[tokfg,line width=.45pt,rounded corners=2.5pt] (0.732,-2.402) rectangle (1.482,-1.952);
\draw[tokmuted,line width=.45pt,rounded corners=3.4pt] (0.629,-2.504) rectangle (1.884,-1.849);
\node[vdot] at (1.733,0.000) {};
\draw[tokfg,line width=.4pt,double,double distance=1.1pt] (2.213,0.000) -- (2.513,0.000);
\draw[tokfg,line width=.5pt] (2.513,0.000) -- (2.813,0.000);
\fill[tokfg] (2.213,0.000) circle (0.061);
\fill[tokfg] (2.513,0.000) circle (0.061);
\fill[tokfg] (2.813,0.000) circle (0.061);
\draw[tokfg,line width=.45pt,rounded corners=1.6pt] (2.391,-0.123) rectangle (2.636,0.123);
\draw[tokfg,line width=.45pt,rounded corners=2.5pt] (2.288,-0.225) rectangle (3.038,0.225);
\draw[tokmuted,line width=.45pt,rounded corners=3.4pt] (1.886,-0.327) rectangle (3.141,0.327);
\node[vdot] at (0.867,1.501) {};
\draw[tokfg,line width=.4pt,double,double distance=1.1pt] (0.957,2.177) -- (1.257,2.177);
\draw[tokfg,line width=.5pt] (1.257,2.177) -- (1.557,2.177);
\fill[tokfg] (0.957,2.177) circle (0.061);
\fill[tokfg] (1.257,2.177) circle (0.061);
\fill[tokfg] (1.557,2.177) circle (0.061);
\draw[tokfg,line width=.45pt,rounded corners=1.6pt] (1.434,2.054) rectangle (1.679,2.299);
\draw[tokfg,line width=.45pt,rounded corners=2.5pt] (1.032,1.952) rectangle (1.782,2.402);
\draw[tokmuted,line width=.45pt,rounded corners=3.4pt] (0.629,1.849) rectangle (1.884,2.504);
\end{tikzpicture}

%% file: fig/assoc-H4.tex
\definecolor{tokfg}{HTML}{1B1F27}
\definecolor{tokmuted}{HTML}{5B6474}
\definecolor{tokline}{HTML}{CFD4DC}
\definecolor{tokfaint}{HTML}{E3E6EC}
\definecolor{toksoft}{HTML}{DDE3EC}
\definecolor{tokbx}{HTML}{2A6FC9}
\definecolor{tokby}{HTML}{D4582A}
\definecolor{toktw}{HTML}{12805F}
\definecolor{toktwsoft}{HTML}{D7EFE6}

\tikzset{
  kedge/.style={draw=tokline,line width=1.71pt,line cap=round},
  kedgetw/.style={draw=toktw,opacity=.5,line width=1.71pt,line cap=round},
  vdot/.style={circle,fill=tokfg,draw=white,line width=.5pt,inner sep=0pt,minimum size=0.360cm},
  tdot/.style={circle,fill=toktw,draw=white,line width=.5pt,inner sep=0pt,minimum size=0.258cm},
  ttdot/.style={rectangle,rounded corners=1pt,fill=toktw,draw=white,line width=.5pt,inner sep=0pt,minimum size=0.270cm},
  sq/.style={fill=toktwsoft,draw=none},
  grid/.style={draw=toktw,opacity=.25,line width=.4pt},
  lab/.style={font=\small,text=tokfg},
  mlab/.style={font=\footnotesize,text=tokmuted},
}

\begin{tikzpicture}
\fill[sq] (2.316,-0.031) -- (1.109,1.215) -- (-2.316,-2.872) -- (-1.109,-4.118) -- cycle;
\draw[grid] (1.914,0.385) -- (-1.511,-3.703);
\draw[grid] (1.511,0.800) -- (-1.914,-3.288);
\draw[grid] (1.888,-0.542) -- (0.681,0.704);
\draw[grid] (1.460,-1.053) -- (0.252,0.193);
\draw[grid] (1.032,-1.564) -- (-0.176,-0.317);
\draw[grid] (0.604,-2.075) -- (-0.604,-0.828);
\draw[grid] (0.176,-2.586) -- (-1.032,-1.339);
\draw[grid] (-0.252,-3.096) -- (-1.460,-1.850);
\draw[grid] (-0.681,-3.607) -- (-1.888,-2.361);
\draw[draw=tokline,opacity=0.59,line width=1.71pt,line cap=round] (2.448,-2.345) -- (1.307,-3.708);
\draw[draw=toktw,opacity=0.31,line width=1.71pt,line cap=round] (2.448,-2.345) -- (-1.175,1.393);
\draw[draw=toktw,opacity=0.32,line width=1.71pt,line cap=round] (2.448,-2.345) -- (3.524,0.175);
\draw[draw=toktw,opacity=0.32,line width=1.71pt,line cap=round] (1.307,-3.708) -- (-2.316,0.031);
\draw[draw=tokline,opacity=0.68,line width=1.71pt,line cap=round] (-1.175,1.393) -- (-2.316,0.031);
\draw[draw=tokline,opacity=0.70,line width=1.71pt,line cap=round] (-0.637,2.653) -- (-1.175,1.393);
\draw[draw=toktw,opacity=0.36,line width=1.71pt,line cap=round] (3.524,0.175) -- (1.109,2.667);
\draw[draw=toktw,opacity=0.37,line width=1.71pt,line cap=round] (1.109,2.667) -- (-0.637,2.653);
\draw[draw=toktw,opacity=0.38,line width=1.71pt,line cap=round] (-1.109,-4.118) -- (1.307,-3.708);
\draw[draw=tokline,opacity=0.77,line width=1.71pt,line cap=round] (-1.241,2.550) -- (-0.637,2.653);
\draw[draw=tokline,opacity=0.77,line width=1.71pt,line cap=round] (-3.524,-0.175) -- (-2.316,0.031);
\draw[draw=toktw,opacity=0.39,line width=1.71pt,line cap=round] (3.524,0.175) -- (2.316,-0.031);
\draw[draw=tokline,opacity=0.80,line width=1.71pt,line cap=round] (0.505,2.564) -- (1.109,2.667);
\draw[draw=toktw,opacity=0.41,line width=1.71pt,line cap=round] (0.505,2.564) -- (-1.241,2.550);
\draw[draw=tokline,opacity=0.83,line width=1.71pt,line cap=round] (-3.524,-0.175) -- (-1.241,2.550);
\draw[draw=toktw,opacity=0.43,line width=1.71pt,line cap=round] (1.109,1.215) -- (0.505,2.564);
\draw[draw=toktw,opacity=0.44,line width=1.71pt,line cap=round] (2.316,-0.031) -- (1.109,1.215);
\draw[draw=toktw,opacity=0.45,line width=1.71pt,line cap=round] (-1.109,-4.118) -- (2.316,-0.031);
\draw[draw=toktw,opacity=0.45,line width=1.71pt,line cap=round] (-2.316,-2.872) -- (-3.524,-0.175);
\draw[draw=toktw,opacity=0.46,line width=1.71pt,line cap=round] (-2.316,-2.872) -- (1.109,1.215);
\draw[draw=toktw,opacity=0.47,line width=1.71pt,line cap=round] (-1.109,-4.118) -- (-2.316,-2.872);
\node[vdot,opacity=0.58] at (2.448,-2.345) {};
\node[tdot,opacity=0.59] at (2.086,-1.971) {};
\node[tdot,opacity=0.59] at (2.556,-2.093) {};
\node[tdot,opacity=0.60] at (1.724,-1.597) {};
\node[tdot,opacity=0.60] at (2.663,-1.841) {};
\node[vdot,opacity=0.60] at (1.307,-3.708) {};
\node[tdot,opacity=0.61] at (1.361,-1.224) {};
\node[tdot,opacity=0.61] at (0.944,-3.334) {};
\node[tdot,opacity=0.61] at (2.771,-1.589) {};
\node[tdot,opacity=0.61] at (0.999,-0.850) {};
\node[tdot,opacity=0.62] at (0.582,-2.960) {};
\node[tdot,opacity=0.62] at (0.637,-0.476) {};
\node[tdot,opacity=0.63] at (2.878,-1.337) {};
\node[tdot,opacity=0.63] at (0.220,-2.586) {};
\node[tdot,opacity=0.63] at (0.275,-0.102) {};
\node[tdot,opacity=0.64] at (1.065,-3.749) {};
\node[tdot,opacity=0.64] at (2.986,-1.085) {};
\node[tdot,opacity=0.64] at (-0.143,-2.212) {};
\node[tdot,opacity=0.64] at (-0.088,0.272) {};
\node[tdot,opacity=0.65] at (-0.505,-1.838) {};
\node[tdot,opacity=0.65] at (3.094,-0.833) {};
\node[tdot,opacity=0.65] at (-0.450,0.646) {};
\node[tdot,opacity=0.66] at (-0.867,-1.465) {};
\node[tdot,opacity=0.66] at (-0.812,1.019) {};
\node[tdot,opacity=0.66] at (3.201,-0.581) {};
\node[tdot,opacity=0.66] at (-1.229,-1.091) {};
\node[tdot,opacity=0.67] at (0.824,-3.790) {};
\node[vdot,opacity=0.67] at (-1.175,1.393) {};
\node[tdot,opacity=0.67] at (3.309,-0.329) {};
\node[tdot,opacity=0.67] at (-1.592,-0.717) {};
\node[tdot,opacity=0.68] at (-1.954,-0.343) {};
\node[tdot,opacity=0.68] at (3.416,-0.077) {};
\node[vdot,opacity=0.69] at (-2.316,0.031) {};
\node[vdot,opacity=0.70] at (3.524,0.175) {};
\node[tdot,opacity=0.70] at (0.582,-3.831) {};
\node[tdot,opacity=0.71] at (2.920,0.798) {};
\node[tdot,opacity=0.72] at (2.316,1.421) {};
\node[vdot,opacity=0.73] at (-0.637,2.653) {};
\node[tdot,opacity=0.73] at (-0.462,2.655) {};
\node[tdot,opacity=0.73] at (-0.288,2.656) {};
\node[tdot,opacity=0.73] at (0.341,-3.872) {};
\node[tdot,opacity=0.73] at (-0.113,2.657) {};
\node[tdot,opacity=0.74] at (3.222,0.123) {};
\node[tdot,opacity=0.74] at (0.061,2.659) {};
\node[tdot,opacity=0.74] at (1.712,2.044) {};
\node[tdot,opacity=0.74] at (0.236,2.660) {};
\node[tdot,opacity=0.74] at (0.410,2.661) {};
\node[tdot,opacity=0.75] at (0.585,2.663) {};
\node[tdot,opacity=0.75] at (0.760,2.664) {};
\node[tdot,opacity=0.75] at (0.934,2.666) {};
\node[vdot,opacity=0.75] at (1.109,2.667) {};
\node[tdot,opacity=0.77] at (0.099,-3.913) {};
\node[tdot,opacity=0.78] at (2.920,0.072) {};
\node[tdot,opacity=0.80] at (-0.143,-3.954) {};
\node[vdot,opacity=0.81] at (-1.241,2.550) {};
\node[tdot,opacity=0.81] at (-1.066,2.552) {};
\node[tdot,opacity=0.81] at (-0.892,2.553) {};
\node[tdot,opacity=0.82] at (-0.717,2.555) {};
\node[tdot,opacity=0.82] at (2.618,0.021) {};
\node[tdot,opacity=0.82] at (-0.542,2.556) {};
\node[tdot,opacity=0.82] at (-0.368,2.557) {};
\node[tdot,opacity=0.82] at (-0.193,2.559) {};
\node[tdot,opacity=0.83] at (-0.019,2.560) {};
\node[tdot,opacity=0.83] at (0.156,2.561) {};
\node[tdot,opacity=0.83] at (-0.384,-3.995) {};
\node[tdot,opacity=0.83] at (0.330,2.563) {};
\node[vdot,opacity=0.84] at (0.505,2.564) {};
\node[tdot,opacity=0.85] at (0.656,2.227) {};
\node[vdot,opacity=0.85] at (-3.524,-0.175) {};
\node[vdot,opacity=0.86] at (2.316,-0.031) {};
\node[tdot,opacity=0.86] at (0.807,1.890) {};
\node[tdot,opacity=0.86] at (-0.626,-4.036) {};
\node[tdot,opacity=0.87] at (-3.403,-0.444) {};
\node[tdot,opacity=0.87] at (1.888,-0.542) {};
\node[tdot,opacity=0.87] at (1.914,0.385) {};
\node[tdot,opacity=0.88] at (0.958,1.553) {};
\node[tdot,opacity=0.88] at (-3.282,-0.714) {};
\node[tdot,opacity=0.88] at (1.460,-1.053) {};
\node[ttdot,opacity=0.88] at (1.486,-0.126) {};
\node[tdot,opacity=0.88] at (1.511,0.800) {};
\node[tdot,opacity=0.88] at (1.032,-1.564) {};
\node[ttdot,opacity=0.89] at (1.057,-0.637) {};
\node[tdot,opacity=0.89] at (-3.162,-0.984) {};
\node[ttdot,opacity=0.89] at (1.083,0.289) {};
\node[vdot,opacity=0.89] at (1.109,1.215) {};
\node[tdot,opacity=0.89] at (0.604,-2.075) {};
\node[ttdot,opacity=0.89] at (0.629,-1.148) {};
\node[ttdot,opacity=0.90] at (0.655,-0.222) {};
\node[tdot,opacity=0.90] at (-0.867,-4.077) {};
\node[tdot,opacity=0.90] at (-3.041,-1.254) {};
\node[tdot,opacity=0.90] at (0.681,0.704) {};
\node[tdot,opacity=0.90] at (0.176,-2.586) {};
\node[ttdot,opacity=0.90] at (0.201,-1.659) {};
\node[ttdot,opacity=0.90] at (0.227,-0.733) {};
\node[tdot,opacity=0.91] at (0.252,0.193) {};
\node[tdot,opacity=0.91] at (-2.920,-1.524) {};
\node[tdot,opacity=0.91] at (-0.252,-3.096) {};
\node[ttdot,opacity=0.91] at (-0.227,-2.170) {};
\node[ttdot,opacity=0.91] at (-0.201,-1.244) {};
\node[tdot,opacity=0.91] at (-0.176,-0.317) {};
\node[tdot,opacity=0.92] at (-2.799,-1.793) {};
\node[tdot,opacity=0.92] at (-0.681,-3.607) {};
\node[ttdot,opacity=0.92] at (-0.655,-2.681) {};
\node[ttdot,opacity=0.92] at (-0.629,-1.755) {};
\node[tdot,opacity=0.92] at (-0.604,-0.828) {};
\node[tdot,opacity=0.93] at (-2.679,-2.063) {};
\node[vdot,opacity=0.93] at (-1.109,-4.118) {};
\node[ttdot,opacity=0.93] at (-1.083,-3.192) {};
\node[ttdot,opacity=0.93] at (-1.057,-2.266) {};
\node[tdot,opacity=0.93] at (-1.032,-1.339) {};
\node[tdot,opacity=0.94] at (-2.558,-2.333) {};
\node[tdot,opacity=0.94] at (-1.511,-3.703) {};
\node[ttdot,opacity=0.94] at (-1.486,-2.777) {};
\node[tdot,opacity=0.94] at (-1.460,-1.850) {};
\node[tdot,opacity=0.95] at (-2.437,-2.603) {};
\node[tdot,opacity=0.95] at (-1.914,-3.288) {};
\node[tdot,opacity=0.95] at (-1.888,-2.361) {};
\node[vdot,opacity=0.96] at (-2.316,-2.872) {};
\node[lab,anchor=north east,inner sep=3pt,fill=white,fill opacity=.8,text opacity=1,rounded corners=1pt] at (2.316,-0.031) {$\{1,12,123\}$};
\node[lab,anchor=south west,inner sep=3pt,fill=white,fill opacity=.8,text opacity=1,rounded corners=1pt] at (-1.109,-4.118) {$\{1,4,12\}$};
\node[lab,anchor=north east,inner sep=3pt,fill=white,fill opacity=.8,text opacity=1,rounded corners=1pt] at (1.109,1.215) {$\{2,12,123\}$};
\node[lab,anchor=south west,inner sep=3pt,fill=white,fill opacity=.8,text opacity=1,rounded corners=1pt] at (-2.316,-2.872) {$\{2,4,12\}$};
\node[vdot] at (-3.524,-4.818) {}; \node[mlab,anchor=west] at (-3.344,-4.818) {nested cube (vertex)};
\node[tdot] at (1.076,-4.818) {}; \node[mlab,anchor=west] at (1.256,-4.818) {non-nested cube on a segment};
\node[ttdot] at (-3.524,-5.268) {}; \node[mlab,anchor=west] at (-3.344,-5.268) {non-nested cube in the interior of a square};
\end{tikzpicture}

%% file: fig/square-H4.tex
\definecolor{tokfg}{HTML}{1B1F27}
\definecolor{tokmuted}{HTML}{5B6474}
\definecolor{tokline}{HTML}{CFD4DC}
\definecolor{tokfaint}{HTML}{E3E6EC}
\definecolor{toksoft}{HTML}{DDE3EC}
\definecolor{tokbx}{HTML}{2A6FC9}
\definecolor{tokby}{HTML}{D4582A}
\definecolor{toktw}{HTML}{12805F}
\definecolor{toktwsoft}{HTML}{D7EFE6}

\tikzset{
  kedge/.style={draw=tokline,line width=0.92pt,line cap=round},
  kedgetw/.style={draw=toktw,opacity=.5,line width=0.92pt,line cap=round},
  vdot/.style={circle,fill=tokfg,draw=white,line width=.5pt,inner sep=0pt,minimum size=0.195cm},
  tdot/.style={circle,fill=toktw,draw=white,line width=.5pt,inner sep=0pt,minimum size=0.140cm},
  ttdot/.style={rectangle,rounded corners=1pt,fill=toktw,draw=white,line width=.5pt,inner sep=0pt,minimum size=0.146cm},
  sq/.style={fill=toktwsoft,draw=none},
  grid/.style={draw=toktw,opacity=.25,line width=.4pt},
  lab/.style={font=\small,text=tokfg},
  mlab/.style={font=\footnotesize,text=tokmuted},
}

\begin{tikzpicture}
\fill[sq] (0,0) rectangle (6.9,3.0);
\draw[grid] (0.863,0) -- (0.863,3.0);
\draw[grid] (1.725,0) -- (1.725,3.0);
\draw[grid] (2.588,0) -- (2.588,3.0);
\draw[grid] (3.450,0) -- (3.450,3.0);
\draw[grid] (4.312,0) -- (4.312,3.0);
\draw[grid] (5.175,0) -- (5.175,3.0);
\draw[grid] (6.038,0) -- (6.038,3.0);
\draw[grid] (0,1.000) -- (6.9,1.000);
\draw[grid] (0,2.000) -- (6.9,2.000);
\draw[kedgetw] (0,0) -- (6.9,0); \draw[kedgetw] (0,3.0) -- (6.9,3.0); \draw[kedgetw] (0,0) -- (0,3.0); \draw[kedgetw] (6.9,0) -- (6.9,3.0);
\node[tdot] at (0.863,0) {}; \node[tdot] at (0.863,3.0) {};
\node[tdot] at (1.725,0) {}; \node[tdot] at (1.725,3.0) {};
\node[tdot] at (2.588,0) {}; \node[tdot] at (2.588,3.0) {};
\node[tdot] at (3.450,0) {}; \node[tdot] at (3.450,3.0) {};
\node[tdot] at (4.312,0) {}; \node[tdot] at (4.312,3.0) {};
\node[tdot] at (5.175,0) {}; \node[tdot] at (5.175,3.0) {};
\node[tdot] at (6.038,0) {}; \node[tdot] at (6.038,3.0) {};
\node[tdot] at (0,1.000) {}; \node[tdot] at (6.9,1.000) {};
\node[tdot] at (0,2.000) {}; \node[tdot] at (6.9,2.000) {};
\node[ttdot] at (0.863,1.000) {};
\node[ttdot] at (0.863,2.000) {};
\node[ttdot] at (1.725,1.000) {};
\node[ttdot] at (1.725,2.000) {};
\node[ttdot] at (2.588,1.000) {};
\node[ttdot] at (2.588,2.000) {};
\node[ttdot] at (3.450,1.000) {};
\node[ttdot] at (3.450,2.000) {};
\node[ttdot] at (4.312,1.000) {};
\node[ttdot] at (4.312,2.000) {};
\node[ttdot] at (5.175,1.000) {};
\node[ttdot] at (5.175,2.000) {};
\node[ttdot] at (6.038,1.000) {};
\node[ttdot] at (6.038,2.000) {};
\node[vdot] at (0.000,0.000) {};
\node[lab,anchor=east,inner sep=5pt] at (0.000,0.000) {$\{1,12,123\}$};
\node[vdot] at (6.900,0.000) {};
\node[lab,anchor=west,inner sep=5pt] at (6.900,0.000) {$\{1,4,12\}$};
\node[vdot] at (0.000,3.000) {};
\node[lab,anchor=east,inner sep=5pt] at (0.000,3.000) {$\{2,12,123\}$};
\node[vdot] at (6.900,3.000) {};
\node[lab,anchor=west,inner sep=5pt] at (6.900,3.000) {$\{2,4,12\}$};
\node[mlab,anchor=north] at (3.450,-0.85) {segment $\{S,12,1\}$};
\node[mlab,anchor=south] at (3.450,3.320) {segment $\{S,12,2\}$};
\node[mlab,anchor=east,align=right] at (-1.9,1.500) {segment\\$\{S,12,123\}$};
\node[mlab,anchor=west,align=left] at (7.200,1.500) {segment\\$\{S,12,4\}$};
\node[mlab] at (0.863,-0.42) {$2$};
\node[mlab] at (1.725,-0.42) {$3$};
\node[mlab] at (2.588,-0.42) {$4$};
\node[mlab] at (3.450,-0.42) {$5$};
\node[mlab] at (4.312,-0.42) {$6$};
\node[mlab] at (5.175,-0.42) {$7$};
\node[mlab] at (6.038,-0.42) {$8$};
\node[mlab,anchor=east] at (-0.1,-0.42) {$k_S=$};
\node[mlab,anchor=east] at (-0.3,1.000) {$k_{12}=2$};
\node[mlab,anchor=east] at (-0.3,2.000) {$k_{12}=3$};
\end{tikzpicture}